\documentclass[reqno,11pt,oneside]{amsart}

\usepackage{url,amsmath,amsthm,amscd,graphicx,psfrag,amssymb,hyperref,setspace}
\usepackage[all]{xy}

\def\tX{{\tilde X}}

\def\MWt{MW_{\rm tors}}
\def\MWl{MW_{\rm lat}}
\def\Aut{{\rm Aut}}

\usepackage[scale=0.73]{geometry}
\let\a=\alpha   \let\b=\beta       \let\d=\delta
    \let\h=\eta          
   \let\l=\lambda  \let\m=\mu      
\let\n=\nu      \let\x=\xi      \let\p=\pi      \let\r=\rho     \let\s=\sigma
\let\t=\tau     \let\o=\omega   \let\c=\chi         
    \let\PH=\Phi        \let\O=\Omega   
      \let\TH=\Theta

 \def\cb{{\mathcal B}}             \def\co{{\mathcal O}} \def\cp{{\mathcal P}}  \def\car{{\mathcal R}} \def\cs{{\mathcal S}}       

 \def\IC{{\mathbb C}} \def\IP{{\mathbb P}}
  
\def\IZ{{\mathbb Z}} \def\IQ{{\mathbb Q}}

\newtheorem{theorem}{Theorem}[section]
\newtheorem{lemma}[theorem]{Lemma}
\newtheorem{corollary}[theorem]{Corollary}
\newtheorem{proposition}[theorem]{Proposition}

\theoremstyle{definition}
\newtheorem{definition}[theorem]{Definition}

\newtheorem*{acknowledgments}{Acknowledgments}

\theoremstyle{remark}
\newtheorem{remark}[theorem]{Remark}

\def\plb#1 #2 {Phys. Lett. {\bf B#1} #2 }
\def\phr#1 #2 {Phys. Rep. {\bf  #1} #2 }
\def\npb#1 #2 {Nucl. Phys. {\bf B#1} #2 }
\def\aph#1 #2 {Ann. Phys. {\bf #1} #2 }
\def\jmp#1 #2 {J. Math. Phys. {\bf #1} #2 }
\def\jgp#1 #2 {J. Geom. Phys. {\bf #1} #2 }
\def\prd#1 #2 {Phys. Rev. {\bf D#1} #2 }
\def\prl#1 #2 {Phys. Rev. Lett. {\bf #1} #2 }
\def\rmp#1 #2 {Rev. Mod. Phys.  {\bf #1} #2 }
\def\zpc#1 {Z. Phys. {\bf #1C} }
\def\cmp#1 #2 {Commun. Math. Phys. {\bf #1} #2 }
\def\cqg#1 #2 {Class.Quant.Grav. {\bf #1} #2 }
\def\mpl#1 {Mod. Phys. Lett. {\bf A#1} }
\def\cpc#1 {Computer Phys. Commun. {\bf #1} }
\def\ijmp#1 {Int. J. Mod. Phys. {\bf A#1} }
\def\ijmpC#1 {Int. J. Mod. Phys. {\bf C#1} }
\def\atmp#1 {Adv. Theor. Math. Phys. {\bf #1} }

\numberwithin{equation}{section}

\let\0=\over
\def\2{{1\over2}}
\let\6=\partial
\def\({\left(}       \def\){\right)}

\let\bra=\langle        \let\ket=\rangle        \def\<#1\>{\bra #1 \ket}
   \def\Im{\text{Im}}

\begin{document}

\title{On a class of non-simply connected Calabi-Yau threefolds}

\author{Vincent Bouchard}
\author{Ron Donagi}

\address{
Perimeter Institute\\
31 Caroline St. N.\\
Waterloo, Ontario\\
N2L 2Y5
}

\email{vbouchard@perimeterinstitute.ca}

\address{University of Pennsylvania\\
Department of Mathematics\\
David Rittenhouse Lab.\\
209 South 33rd Street\\
Philadelphia, PA 19104-6395
}

\email{donagi@math.upenn.edu}

\begin{abstract}
We obtain a detailed classification for a class of non-simply connected Calabi-Yau threefolds which are of
potential interest for a wide range of problems in string
phenomenology. These threefolds arise as quotients of Schoen's
Calabi-Yau threefolds, which are fiber products over $\IP^1$ of two
rational elliptic surfaces. The quotient is by a freely acting finite
abelian group preserving the fibrations. Our work involves a classification of restricted finite automorphism groups of rational elliptic surfaces.
\end{abstract}

\maketitle

\makeatletter
\renewcommand{\l@section}{\@tocline{1}{2pt}{1pc}{}{}}
\renewcommand{\tocsection}[3]{%
  \indentlabel{\bf \@ifnotempty{#2}{ \ignorespaces#1 #2.\quad}} \bf #3}
\renewcommand{\l@subsection}{\@tocline{2}{0pt}{4pc}{5pc}{}}
\makeatother

\tableofcontents

\section{Introduction}

The goal of this paper is to construct and study a class of non-simply connected 
Calabi-Yau threefolds which are of potential interest for a wide range of problems in string phenomenology. These threefolds arise as quotients of Schoen's Calabi-Yau threefolds, 
which are fiber products over $\IP^1$ of two rational elliptic surfaces 
$\beta:B \to {\IP^1},  \beta':B' \to {\IP^1}$.  
Our goal is to classify pairs of actions $\phi_{B}, \phi_{B'}$ 
of a finite abelian group $G$ on the two rational elliptic surfaces 
$\beta:B \to {\IP^1},  \beta':B' \to {\IP^1}$, compatible with the same action 
$\phi_{\IP^1}$ of $G$ on $\IP^1$, such that the resulting action 
$\phi_{\tilde X} := \phi_B \times_{\phi_{\IP^1}}  \phi_{B'}$
of $G$ on the fiber product 
${\tilde X} := B \times_{\IP^1} B'= \{ (p,p') \in B \times B' | \b'(p') = \b (p) \} $
is free, yielding a non-simply connected quotient Calabi-Yau threefold $X:= \tX / G$ 
with fundamental group $\pi_1(X) = G$. This goal is accomplished in theorem \ref{t:main}.

The direct motivation for this came from the problem of deriving the Standard Model 
of particle physics as the low energy limit of an heterotic string compactification.
We discuss some of the issues involved and review the known results in section 1.1. 
We then describe the precise version of our problem in section 1.2. 
The remainder of the paper is outlined in section 1.3.

\subsection{Physics motivation}

The main motivation for this work comes from physics. An important open problem in string theory is to find a low-energy limit of the theory that reproduces observed physics, such as the Standard Model of particle physics. There are many approaches to this problem; one of which consists in studying four-dimensional compactifications of the so-called $E_8 \times E_8$ heterotic string theory.

Roughly speaking, an $E_8 \times E_8$ heterotic vacuum is characterized by a triple $(X,V,G)$, where $X$ is a Calabi-Yau threefold and $V$ is a stable vector bundle on $X$ with structure group $G \subseteq E_8$.\footnote{For the readers familiar with this construction, here we implicitly consider the hidden bundle in the second $E_8$ factor to be trivial.} In this model, the low-energy physics translates into topological and geometrical properties of this triple. It turns out that obtaining realistic physics out of $E_8 \times E_8$ compactifications on smooth simply connected Calabi-Yau threefolds $X$ runs into various difficulties.

Many problems may be circumvented by considering instead smooth non-simply connected Calabi-Yau threefolds. For instance, an heterotic vacuum on a smooth Calabi-Yau threefold $X$ with fundamental group $\p_1(X) = \IZ_2$ reproducing precisely the massless spectrum of the Minimal Supersymmetric Standard Model of particle physics has recently been found \cite{Bouchard:2005ag, Bouchard:2006dn}. It is plausible that many more such realistic vacua exist in the landscape of string vacua.

Hence, the holy grail of connecting string theory to the real world motivates the search for non-simply connected Calabi-Yau threefolds. Sporadic examples of non-simply connected Calabi-Yau threefolds already exist in the literature, but they are rather isolated compared to the wild jungle of simply connected Calabi-Yau threefolds. In this paper we undertake the study of a large class of non-simply connected Calabi-Yau threefolds, based on Schoen's construction of smooth simply connected Calabi-Yau threefolds as (resolutions of) fiber products of two rational elliptic surfaces \cite{Schoen:1988}.

\subsection{Geometric setup}

Our aim is to construct smooth Calabi-Yau threefolds $X$ with fundamental group $\p_1(X) = G$, where $G$ is a finite abelian group. Every such $X$ has a smooth simply connected universal cover $\tX$, which comes with a freely acting group of automorphisms $G_{\tX} \cong G$. Moreover, by pulling back the Calabi-Yau metric of $X$ we see that the universal cover $\tX$ is also Calabi-Yau.

Conversely, let $\tX$ be a smooth simply connected Calabi-Yau threefold, with a freely acting 
finite abelian group of automorphisms $G_{\tX}$. It turns out that $G_{\tX}$ automatically preserves 
the volume form\footnote{
The argument, due to Beauville, goes as follows. Let $\tX$ be a Calabi--Yau threefold, by which we mean a threefold with $SU(3)$ holonomy; then $h^{1,0} (\tX)=h^{2,0} (\tX) = \chi(\co_{\tX}) = 0$. If $G_{\tX}$ acts freely on $X$, then the quotient manifold $X = \tX / G_{\tX}$ also satisfies $h^{1,0} (X)=h^{2,0} (X) = \chi(\co_{X}) = 0$. It follows that $h^{3,0}(X) = 1$, hence there is a nonzero $(3,0)$-form $\O$ on $X$. Since its pullback on $\tX$ is nowhere vanishing, $\O$ must also be nowhere vanishing on $X$; therefore $X$ is Calabi--Yau \cite{Beauville}.}, so the quotient manifold $X = \tX / G_{\tX}$ is a smooth Calabi-Yau threefold with 
fundamental group $\p_1(X) \cong G_{\tX}$.

In this paper we use the second point of view to construct a large class of non-simply connected 
Calabi-Yau threefolds. Namely, we classify all Calabi-Yau threefolds $\tX$ constructed as a smooth 
fiber product of two rational elliptic surfaces and admitting a freely acting finite abelian group 
of automorphisms $G_{\tX}$. 
Since we will eventually also need to use the spectral construction to build Standard Model bundles 
on $X$, we require that $G_{\tX}$ preserves the elliptic fibration of $\tX$, so that the 
quotient threefold $X$ is torus-fibered.\footnote{By torus-fibered we mean that $X$ has a fibration with 
generic fiber an elliptic curve, but the fibration does not necessarily have a section.} 

\vspace{1em}

Let us now describe in more details the family of Calabi-Yau threefolds we are interested in. Let $B$ and $B'$ be rational elliptic surfaces. Consider the fiber product $\tX := B \times_{\IP^1} B'$. That is, construct the threefold
\begin{equation}
{\tilde X} = \{ (p,p') \in B \times B' | \b'(p') = \b (p) \},
\end{equation}
 where $\b: B \to \IP^1$ and $\b': B' \to \IP^1$ are the elliptic fibrations of the rational elliptic surfaces $B$ and $B'$. $\tilde X$ can also be described by the commutative diagram
\begin{equation}\label{e:fibration}
\xymatrix{
& {\tilde X} \ar[dl]_{\pi'} \ar[dr]^{\pi} \\
B \ar[dr]^{\b} && B' \ar[dl]_{\b'} \\
& \IP^1
}
\end{equation}

The rational elliptic surfaces $B,B'$ are smooth, but some of the fibers of $B,B'$ over $\IP^1$ must be singular.
Denote by $S$ (resp. $S'$) the set of images of the singular fibers of $B$ (resp. $B'$) in $\IP^1$. 
The fiber product $\tX$ will be smooth except at points $(p,p')$ 
where both $p$ and $p'$ are singular points of their respective fibers.
The image in $\IP^1$ of the singularities is therefore the intersection $S'' = S \cap S'$.
In other words, $\tX$ is smooth if and only if $S'' = \emptyset$, that is, no singular fibers 
of $B$ and $B'$ are paired in the fiber product.

It was shown by Schoen \cite{Schoen:1988} that the smooth $\tX$ are Calabi-Yau. 
In fact, Schoen considered the more general case when $\tX$ has ordinary double point singularities. 
More precisely, he studied the case where for any $p \in S''$, the singular fibers 
$\b^{-1}(p)$ and $(\b')^{-1}(p)$ are semistable --- of type $I_k$, $k \geq 0$, in Kodaira's notation. Then, if the singular fibers above the points in $S''$ satisfy certain additional requirements, the minimal resolution of the singular threefold $\tX$ is also Calabi-Yau. However, in this paper we will not consider these additional smooth Calabi-Yau threefolds and focus on smooth fiber products.

Since the Euler characteristic of $T^2$ is $0$, it is clear that the Euler characteristic 
of the fiber product is given by the sum of the Euler characteristics of the singular fibers 
above the points in $S''$. In particular, if $S'' = \emptyset$, then the fiber product $\tX$ is smooth, 
and $\c (\tilde X) = 0$. One can show that $\tX$ is simply connected and its Hodge diamond is given by
\begin{equation}
\begin{matrix}
&&&1&&&\\
&&0&&0&&\\
&0&&19&&0&\\
1&&19&&19&&1\\
&0&&19&&0&\\
&&0&&0&&\\
&&&1&&&
\end{matrix}
\end{equation}
The cohomology group $H^2 (\tilde X, \IZ) \simeq {\rm Pic} (\tilde X)$ of $\tilde X$ is given by
\begin{equation}
H^2 (\tilde X, \IZ) \cong \frac{H^2(B, \IZ) \oplus H^2(B', \IZ) }{H^2( \IP^1, \IZ)}.
\end{equation}

So our aim is to give a complete list of pairs $(\tX,G_{\tX})$ with $\tX$ a smooth fiber product of two rational elliptic surfaces and $G_{\tX}$ a freely acting finite abelian group of automorphisms
preserving the elliptic fibration. Because of the fiber product nature of $\tX$ described by the commutative diagram \eqref{e:fibration}, our strategy will be to reduce the classification of pairs $(\tX,G_{\tX})$ to a classification of finite automorphism groups of rational elliptic surfaces $B$.

\vspace{1em}

An automorphism $\t_\tX : \tX \to \tX$ of a Calabi-Yau threefold $\tX = B \times_{\IP^1} B'$ 
that preserves the elliptic fibration has the form $\t_{\tX} = \t_B \times_{\IP^1} \t_{B'}$, 
where $\t_B$ and $\t_{B'}$ are automorphisms of the rational elliptic surfaces $B$ and $B'$. 
Moreover, since there is no free automorphism of a rational elliptic surface,\footnote{See section \ref{s:fixed} for a short proof of this fact.} it is clear that if $\t_{\tX}$ acts freely on $\tX$, then both $\t_B$ and $\t_{B'}$ 
must act non-trivially on $B$ and $B'$ respectively. Hence, if $\tX$ is a smooth fiber product 
with a freely acting, fibration-preserving automorphism group $G_{\tX}$, then it must be 
the fiber product of two rational elliptic surfaces $B$ and $B'$ with automorphism groups 
$G_B \cong G_{B'} \cong G_{\tX}$.\footnote{Note here that $G_B$ and $G_{B'}$ are not necessarily the full automorphism groups of the rational elliptic surfaces, but may be subgroups thereof.}

The direct consequence of this observation is that the complete list of pairs $(\tX,G_{\tX})$ can be easily reconstructed from the list of pairs $(B,G_B)$. Namely, given the complete list of pairs $(B,G_B)$, one has to consider ``pairs of pairs" $(B,G_B)$ and $(B',G_{B'})$ with $G_B \cong G_{B'}$ and consistent induced actions on the $\IP^1$ base to construct all possible $(\tX,G_{\tX})$. Henceforth, we will be interested in finite automorphism groups $G_B$ of rational elliptic surfaces $B$.

Notice however that not all automorphisms of rational elliptic surfaces $B$ can be lifted to free automorphisms of $\tX$. For $\t_{\tX} = \t_B \times_{\IP^1} \t_{B'}$ to be free on $\tX$, it is clear that the fibers of $B$ and $B'$ intersecting the fixed loci of $\t_B$ and $\t_{B'}$ must not be paired in the fiber product. This gives restrictions on allowed automorphisms of rational elliptic surfaces, which we explore in section \ref{s:restrictions}. 
In those cases where the lifting to a free automorphism is possible, we
typically get infinitely many such free autmorphisms, parametrized by a
finite number of cosets of an infinite group. However, only a finite
number of distinct quotients arise:  automorphisms in the same coset are
related via conjugation by a translation automorphism, hence they result
in isomorphic quotients. This is discussed in section \ref{s:summproc}.

\subsection{Outline}

We start by recalling important properties of rational elliptic surfaces in section \ref{s:prop}. 
We study automorphisms of rational elliptic surfaces in section \ref{s:auto}, 
and explore the properties that they must satisfy in order to lift to free automorphisms on smooth fiber products. 
In particular, we split the task of classifying such automorphisms into separate questions about 
the ``linear" part, which preserves the zero-section, 
and the ``translation" part, which has trivial linearization.
Then, in sections \ref{s:linear} and \ref{s:sections} we study what we call 
\emph{automorphisms of the second kind}, which are automorphisms acting non-trivially on the $\IP^1$ base. 
The list of suitable linearized actions is obtained in Table 4 in section \ref{s:linear}, 
and the corresponding translations are classified in section \ref{s:sections}.
In section \ref{s:class}, we combine our analysis of the previous sections to understand 
finite automorphism groups of rational elliptic surfaces that can be lifted to 
free automorphism groups on smooth fiber products. 
We produce in tables \ref{t:listtriv}--\ref{t:lists} a complete list of pairs $(B,G_B)$, where $G_B$ is a finite automorphism group of $B$ that can be lifted to a free finite automorphism group of $\tX$.
The list of threefolds is obtained from all pairs $B,B'$ 
with the same group $G$ and the same action on $\IP^1$; this is accomplished in our main theorem \ref{t:main}, stated in section \ref{s:examples}. We conclude this section by giving some examples of non-simply connected Calabi-Yau threefolds obtained from our results. Finally, we present in section \ref{s:outlook} some open questions and avenues for future research.

\subsection{Hitchiker's guide to the paper}

Since our calculations are sometimes rather detailed, it may help the
reader to have one brief outline of the entire argument.

We want to classify free, fibration-preserving actions of a finite
group $G$ on a smooth Calabi-Yau fiber product $\tX =B
\times_{\IP^1}B'$ and the quotient Calabi-Yaus $X = \tX / G$. Each
group element acts as an automorphism $\tau_{\tX}$ which itself, by
assumption, is a fiber product $\tau_{\tX} = \tau_{B}
\times_{\tau_{\IP^1}}\tau_{B'}$ where $\tau_{B}$ is an automorphism
of $B$, $ \tau_{B'}$ is an automorphism of $B'$, and they both lift
the same automorphism $\tau_{\IP^1}$ of $\IP^1$. 

The automorphism $\tau_{B}$ is ``affine", in the sense that it does not necessarily
send the zero-section $\sigma$ to itself; it has a {\it
linearization} $\alpha_B$ which does send $\sigma$ to itself. The linearization
lifts the same automorphism $\tau_{\IP^1}$ of $\IP^1$, and in fact
differs from $\tau_B$ by a translation (by a section $\xi$ of $B \to
\IP^1$ ): $\tau_B = t_{\xi} \circ \alpha_B$. So we reduce the
determination of $G$ actions to consideration of the possible
automorphisms $\tau_{\tX}$; these are reduced first to analysis of
the possible $\tau_B$, and then to a study of the possible $\alpha_B$
and $\xi$.

In section 3.2 we introduce notation for the orders of the various
finite order automorphisms:  $n:=\text{order}(\tau_B)$,
$m:=\text{order}(\alpha_B)$, $\bar{m}:=\text{order}(\tau_{\IP^1})$. There are
positive integers $d,k$ such that:
$$n=md=\bar{m}kd.$$
We say that an automorphism is of the {\it first kind} if $m=1$, and
of the {\it second kind} if $m>1$. 

A key notion is that of a {\it
suitable pair} $(B,\tau_B)$ (Definition 3.9): $(B,\tau_B)$ is
suitable if $k=1$ (i.e. $m=\bar{m}$) and either:
\begin{itemize}
\item
$m=1$ (first kind), or:
\item
$m>1$ (second kind), and:
\begin{itemize}
\item the subgroup $\langle \tau_B \rangle$ generated by $\tau_B$
acts {\it freely} on the fiber $f_{\infty}$ above one of the fixed
points $\infty \in \IP^1$ of $\tau_{\IP^1}$, and:
\item
if $d>1$ we also require that $f_{\infty}$ be
smooth.
\end{itemize}
\end{itemize}

In section 3.3 we show that an automorphism $\tau_B$, in conjunction with a similar automorphism $\tau_{B'}$, can be lifted to
a free $\tau_{\tX}$ if and only if the pairs $(B,\tau_B)$ and $(B', \tau_{B'})$ are
suitable.

The next reduction uses the notion (definition \ref{d:sigma})
of a {\it suitable $\sigma$-pair}, i.e. a pair $(B,\alpha_B)$ where
$\alpha_B$ is an automorphism sending $\sigma$ to itself and
satisfying:
\begin{itemize}
\item
$m=\bar{m}$;
\item
$f_{\infty}$ is of Kodaira type $I_{mr}$ for some integer $r \geq 0$;
\item
$\alpha_B$ fixes the neutral component of $f_{\infty}$ pointwise.
\end{itemize}

We show in section 3 that if $(B,\tau_B)$ is a suitable pair, then its
linearization $(B,\alpha_B)$ is a suitable $\sigma$-pair. In
particular, in the cases relevant to us, $m=\bar{m}$, so we will use
$m$ and $d$ (or $m$ and $n=md$) as our basic discrete parameters.

Suitable $\sigma$-pairs are classified in proposition 4.8 and table 4.
In that table we also provide basic information about the rational
elliptic surface $\widehat{\cb}$ obtained as the Kodaira model of the
quotient $B / \bra  \alpha_B \ket$. The classification is done under the
general-position assumption that the fiber $f_{\infty}$ is smooth. As
explained in section 4.3.1, this assumption is harmless: it is satisfied at
the generic point of each of the moduli spaces we are concerned with,
so we do not miss any components by imposing it. The classification is
based on the notion of the {\it deficiency} of a singular fiber of the Kodaira model $\widehat{\cb}$ of the quotient $B / \bra \alpha_B \ket$ (definition
4.2). Topological considerations lead to constraints on fiber
deficiencies (lemmas 4.4, 4.6). These imply that $m$ must  be between
1 and 6,  determine the suitability of $\widehat{\cb}$ in terms of its special
fibers $f_0,f_{\infty}$ over the fixed points $0$ and $\infty$ of $\t_{\IP^1}$ (lemma 4.7), and lead to the proof of
proposition 4.8.

Given $\alpha_B$, for each positive integer $i$ let $\cp_i$ be the
endomorphism of the Mordell-Weil group $MW$ induced by the element
$\sum_{j=0}^i \alpha_B^j$ of the group ring $\IZ[G]$ of $G$. We also
set $\Phi_i(\xi) := \langle \cp_i(\xi), \cp_i(\xi) \rangle$ where $\langle \cdot, \cdot \rangle$ is the height
pairing, reviewed in section 2. Note that $\Phi_m(\xi)=0$ if and only if
$\cp_m(\xi)$ is a torsion section of some finite order $d$. Hence, for any section $\xi \in \ker(\PH_m)$, the automorphism $\tau_B = t_\xi \circ \alpha_B$ has order $n=md$. According to lemma 5.1,
the determination of the kernel of $\Phi_m$ amounts to calculating
the perpendicular to the sublattice $MW^{\alpha_B}$ in $MW$. This
rather tedious calculation is carrried out case by case in lemma 5.2. 

The possible orders $d$ of $\cp_m(\xi)$ must divide the order of the
torsion subgroup $MW_{tors}$ (lemma 3.2), and if the fiber $f_0$ is
smooth (where $0 \in \IP^1$ is the other fixed point of
$\tau_{\IP^1}$) then $d=1$. This leaves a fairly short list of
possibilities for $d>1$, worked out in lemma 5.3. 
We find that $d=1$ for $\xi$ in the subgroup $\ker(\PH_m)_{d=1} =
\cp_m^{-1}(0)$ of $\ker(\PH_m) = \cp_m^{-1}(\MWt)$. In each case there is one
other possible value for $d$, which always turns out to be 2 or 3 and
which is obtained for $\xi$ in the union of the non-trivial cosets of
$\ker(\PH_m)_{d=1}$ in $\ker(\PH_m)$.

Given a suitable $\sigma$-pair $(B,\alpha_B)$ of order $m$, the set of
{\it allowed sections} $AS$ is then given by $\ker(\PH_m)$, minus the cosets of $\Im(1-\alpha_B)$ in $\ker(\PH_m)$ containing sections which intersect the fiber at infinity $f_\infty$ at a torsion point of order less than $n = d m$. In
lemma 3.4 and proposition 5.4 we show that $\xi$ is an allowed section
in this sense if and only if the pair $(B,\tau_B = t_{\xi} \circ
\alpha_B)$ is suitable. 

Quite remarkably, it turns out that allowed sections always exist: in
each of the relevant cases there is at least one coset that satisfies the
requirements, i.e. the sections belonging to it meet the fiber
$f_{\infty}$ at a point of order precisely $n$. This is noted in remark
6.3, and the proof is outlined in the Appendix. The result is that each
possibility allowed by our lattice-theoretic calculations actually
corresponds to a non-empty family of rational elliptic surfaces and thus
produces a family of Schoen threefolds with a free action of the
appropriate finite group.

At this point, all the pieces are ready to be assembled. The two
propositions in section 6 tabulate rational elliptic surfaces with
finite automorphism groups of the first and second kind,
respectively.  These in turn are used in the Main Theorem (7.1) to
produce the desired list of non simply connected Calabi-Yau
threefolds.

\begin{acknowledgments}
We would like to thank Chuck Doran, Antonella Grassi, John McKay, Rick Miranda, Tony Pantev, 
Ulf Persson, Matthias Schuett and Noriko Yui for valuable discussions. We would also like to add special thanks to the referee for relevant comments and suggestions.

V.B. thanks the Mathematics Department of University of Pennsylvania 
for hospitality while part of this work was completed. 
R.D. is partially supported by NSF grants DMS 0104354 and DMS 0612992, 
and an NSF Focused Research Grant DMS 0139799 for ``The Geometry of Superstrings". 
V.B. is partially supported by an NSERC postdoctoral fellowship and by Perimeter Institute for Theoretical Physics. 
Research at Perimeter Institute for Theoretical Physics is supported in part by the Government of Canada through NSERC and by the Province of Ontario through MRI.
\end{acknowledgments}


\section{Preliminaries}\label{s:prop}

Let us start by recalling some important properties of rational elliptic surfaces. Note that everything in the following will be defined over the complex numbers.

\subsection{General properties}
Let $B$ be a rational elliptic surface and denote by $\b : B \to \IP^1$ its elliptic fibration over $\IP^1$. Let $\s:\IP^1 \to B$ be the zero section. 

The Euler characteristic of $B$ is $\c (B) = 12$. Since the Euler characteristic of a smooth elliptic fiber is $\c (T^2) = 0$, $B$ must have singular fibers with Euler characteristics adding up to $12$. As is well known, the complete list of possible fibers has been worked out by Kodaira \cite{Kodaira}: it includes two infinite families --- $I_n$ and $I^*_n$ for $n \geq 0$, and six exceptional cases --- $II$, $III$, $IV$, $II^*$, $III^*$, $IV^*$. Moreover, $B$ is simply connected; therefore the single nontrivial cohomology group is given by $H^2(B,\IZ) \cong \IZ^{10}$.

\subsection{Mordell-Weil group}
The {\em Mordell-Weil group} $MW$ of a rational elliptic surface $B$ is the group of rational sections. A well-known theorem by Mordell and Weil states that $MW$ is finitely generated.

Define the subgroup $T \subset H^2(B,\IZ)$ generated by the cohomology classes of the zero section $\sigma$, the generic fiber $f$ and the irreducible components of the singular fibers not intersecting $\sigma$. 
Then the Mordell-Weil group of $B$ can be identified with the quotient:
\begin{equation}\label{e:MW}
MW \cong H^2(B,\IZ) / T.
\end{equation}

\subsection{Mordell-Weil lattice}
Let $\MWt$ be the torsion subgroup of the Mordell-Weil group. Define
$$
\MWl = MW/\MWt.
$$
Following \cite{Shioda:1990, Oguiso:1991}, we can endow the group $\MWl$ with an additional lattice structure --- hence the subscript \emph{lat}. Let $\m,\n \in MW$ be two sections. Define the height pairing of $\m$ and $\n$ by
\begin{equation}\label{e:hp}
\< \m,\n \> = 1 + \m \s + \n \s - \m \n - \sum_{s_i} {\rm contr}_{s_i} (\m,\n),
\end{equation}
where the sum runs over the singular fibers $s_i$ of $B$, 
multiplication denotes intersection numbers of the sections, and 
${\rm contr}_{s_i} (\m,\n)$ is the entry $({\m}_i,{\n}_i)$
of the inverse of the intersection matrix $T_{s_i}$ of the A-D-E lattice associated with the singular fiber $s_i$. 
Here ${\m}_i, {\n}_i$ denote the component of singular fiber $s_i$ intersecting the sections $\m,\n$ respectively.  
In particular, the numbers ${\rm contr}_{s_i} (\m,\n)$ are rational, and zero whenever $\m$ or $\n$ intersects the neutral component of $s_i$. This defines the structure of a positive-definite lattice on $\MWl$.

Define the narrow Mordell-Weil lattice to be the sublattice 
$\MWl^0$ consisting of sections $\m \in \MWl$ passing through the neutral component of every fiber. For $\m,\n \in \MWl^0$, the height pairing takes the simplified form
$$
\< \m,\n \> = 1 + \m \s + \n \s - \m \n,
$$
which is integer valued. It promotes $\MWl^0$ to the status of a positive-definite even integral lattice \cite{Shioda:1990}. Moreover, $\MWl^0 \cong \MWl^*$, where $\MWl^*$ denotes the dual lattice of $\MWl$.

Recall now a few properties that will be useful later. 

\begin{proposition}[\cite{Oguiso:1991}, proposition 3.5]\label{p:torsion}
A section $\h \in MW$ is torsion if and only if $\langle \h, \h \rangle = 0$. Moreover, for rational elliptic surfaces defined over the complex number field, any non-zero torsion section is disjoint from the zero section. This implies (using \eqref{e:hp}) that a section $\h \in MW$, $\h \neq \s$, is torsion if and only if $\sum_{s_i} {\rm contr}_{s_i} (\h,\h) = 2$.
\end{proposition}

\begin{lemma}[\cite{Shioda:1990}, lemma 10.7]
For any rational elliptic surface, the number of the sections $\x$ which are disjoint from the zero section $\s$ is finite and at most 240. Every such $\x$ satisfies $\langle \x, \x \rangle \leq 2$.
\end{lemma}

\begin{theorem}[\cite{Shioda:1990}, theorem 10.8]\label{t:generators}
The Mordell-Weil group $MW$ of a rational elliptic surface is generated by the sections $\x$ which are disjoint from the zero section $\s$, hence by those satisfying $\langle \x, \x \rangle \leq 2$.
\end{theorem}

Finally, the following result on lattices, which is used in \cite{Oguiso:1991} to prove the theorem above, will also be useful.

\begin{lemma}[\cite{Oguiso:1991}, lemma 5.1]
Let $L$ be a root lattice of type $A_n$, $D_n$ or $E_k$, and let $L^*$ be its dual lattice, namely the weight lattice of type $A_n$, $D_n$ or $E_k$. Then $L^*$ has a basis consisting of minimal vectors of minimal norm
$$
\m(A_n^*) = \frac{n}{n+1},~~~~~~~~~\m(D_n^*) = 1,~~~~~~~~~ \m(E_k^*) = \frac{10-k}{9-k} = \frac{4}{3}, \frac{3}{2}, 2 ~~(k=6,7,8),
$$
except in the case $L^* = D_n^* (n > 4)$ where a vector of norm $n/4$ should be added.
\end{lemma}

\subsection{Weierstrass model}

Any rational elliptic surface $B$ can be described as the minimal resolution of a Weierstrass model $W(\co(1),A_4,A_6)$ over $\IP^1$, where $\co(1)$ is the tautological line bundle on $\IP^1$, and the sections $A_4 \in H^0 (\IP^1, \co(4))$ and $A_6 \in H^0 (\IP^1, \co(6))$ can be described by polynomials in $t$ of order $4$ and $6$ respectively. The Weierstrass model $W(\co(1),A_4,A_6)$ is defined by the equation
\begin{equation}\label{e:weierstrass}
y^2 z = x^3 + A_4(t) x z^2 + A_6(t) z^3,
\end{equation}
where $[x:y:z]$ are homogeneous coordinates\footnote{More precisely \eqref{e:weierstrass} defines a divisor in the threefold given by the $\IP^2$ fibration $P=\IP(\co(2) \oplus \co(3) \oplus \co) \to \IP^1$, where $\IP()$ means projectivization. $[x:y:z]$ is then a section of the bundle $P$.} of $\IP^2$. The zero section $\sigma$ is simply $[0:1:0]$.

The map $\r: B \to W$ from $B$ to its Weierstrass model $W$ is given by blowing down all the irreducible components of the singular fibers of the elliptic fibration not meeting the zero section. Thus, the fibers of $W$ are all irreducible, but $W$ is singular (unless $W = B$); all its singularities are rational doublepoints, and $B$ is the minimal desingularization.

The surfaces $B,W$ are fibered over the same base $\IP^1$.
Let $t_i$ be the points in this $\IP^1$ base where the fibers of $\beta$ are singular. 
The points $t_i$ are given by the vanishing of the discriminant of the elliptic fiber
\begin{equation}
\Delta = 4 A_4(t)^3 + 27 A_6(t)^2.
\end{equation}
As is well known, at each of the $t_i$'s the Kodaira type of the singular fiber in the resolved surface $B$ 
can be read off directly from the order of vanishing of $\Delta$, $A_4$ and $A_6$.

\section{Automorphisms of rational elliptic surfaces}\label{s:auto}

Let us now turn to the study of automorphisms of finite order on rational elliptic surfaces.

Let $B$ be a rational elliptic surface, $\beta: B \to \IP^1$ its elliptic fibration, and $\s: \IP^1 \to B$ the zero section.

\subsection{General properties}

\subsubsection{Fixed locus}\label{s:fixed}


Let $\tau_B$ be a non-trivial, fixed point free automorphism of a smooth rational surface B.
If the order $n$ of  $\tau_B$  is prime, all non-trivial elements in the group $\langle \t_B \rangle$ 
generated by $\tau_B$  are fixed point free, so the quotient $B' = B/\langle \t_B \rangle$
is smooth and the quotient map $q : B \to B'$ is a covering map. But by Luroth's theorem 
$B'$ is rational, hence simply connected, so it cannot admit a non trivial cover. Inductively, 
we see that this holds for any $n>1$: the fixed locus of any non trivial automorphism $\tau_B$ 
of a rational elliptic (or just rational) surface $B$ must be non empty.

\subsubsection{Types of automorphisms}

Any automorphism preserves the canonical class. For a rational elliptic surface $B$, the canonical class is $K_B = -f$, where $f$ is the class of a generic elliptic fiber. Thus, any automorphism of $B$ preserves the fiber class, hence the elliptic fibration. In other words, any automorphism $\t_B : B\to B$ maps fibers to fibers.

As a result, an automorphism of a rational elliptic surface can be of one of two
types: either it leaves each fiber stable, or it permutes fibers and descends 
to a non trivial 
automorphism on $\IP^1$. We call the former \emph{automorphisms of the first kind}, and the latter \emph{automorphisms of the second kind}.

\subsubsection{Group homomorphisms}

Let us now introduce a few useful group homomorphisms. Denote by $\Aut(B)$ the automorphism group of $B$, 
by $\Aut(\IP^1)$ the automorphism group of $\IP^1$, and by $\Aut_\s(B)$ the group of automorphisms of $B$ that preserve the zero section $\s$.

A first group homomorphism that we consider is given by $t : MW \to {\rm Aut}(B)$, which sends a section $\x \in MW$ to the automorphism $t_\x : B \to B$ defined by translation by the section $\x$. The automorphisms $t_\x$ are of the first kind, since they leave each fiber stable.

There is a canonical group homomorphism $\rho: \Aut(B) \to \Aut(\IP^1)$ which sends
$$
\rho: \t_B \mapsto \rho(\t_B):= \t_{\IP^1}.
$$ 
Explicitly, $\t_{\IP^1}$ is specified by $\t_{\IP^1} \circ \b = \b \circ \t_B$.
In other words, $\rho$ associates to any automorphism $\t_B : B \to B$ its induced automorphism on the $\IP^1$ base $\t_{\IP^1}: \IP^1 \to \IP^1$. 
An automorphism $\t_B$ is of the first kind if and only if $\t_{\IP^1} = 1$, while it is of the second kind otherwise.
In particular, $\rho(t_{\xi})=1$ for all $\xi \in MW$.

There is another group homomorphism $\lambda: \Aut(B) \to \Aut_{\s}(B)$ acting by
$$
\lambda: \t_B \mapsto \l(\t_B):=\a_B,
$$
where $\a_B \in \Aut_\s(B)$ 
is the \emph{linearization of $\t_B$}. The precise form of the map goes as follows. Let $\x = \t_B(\s)$ be the section defined by applying $\t_B$ to the zero section $\s$. Then define the linearization $\a_B$ of $\t_B$ by $\a_B = t_{-\x} \circ \t_B$. By definition, $\a_B(\s)=\s$, hence $\a_B \in \Aut_\s (B)$. To see that $\l$ is compatible with the group structure, recall that fiber by fiber it is an isomorphism between genus 1 curves which sends the origin to the origin, hence a group homomorphism.

\begin{remark}
In the following, we will always use the symbols $\t_{\IP^1}$ and $\a_B$ for the images of $\t_B$ under $\rho$ and $\l$ respectively.
\end{remark}

\subsection{Automorphisms of finite order}

Let $\t_B \in \Aut(B)$ be an automorphism of finite order. Denote by $n$ the order of $\t_B$, by $m$ the order of $\a_B$, and by $\bar m$ the order of $\t_{\IP^1}$. We have that
$$
n = m d = {\bar m} k d,
$$
for some integers $d,k > 0$.

\begin{lemma} \label{l:k}
$k$ divides the order of the group of complex multiplications
(automorphisms fixing the origin) of the generic elliptic fiber. In
particular, $k \in \{1,2,3,4,6\}$. $d$ divides the order of the torsion
subgroup $\MWt$.
\end{lemma}

\begin{proof}
Since $\t_{\IP^1} = \rho(\alpha_{\beta})$ has order $\bar m$, 
we get that $(\a_B)^{\bar m}$ must leave each fiber stable. By definition, it must also fix the zero section. Hence it must act on each elliptic fiber separately and fix the zero point; that is, it must act by complex multiplication on the generic elliptic fiber. Consequently, $k = m / {\bar m}$ must divide the order of the group of complex multiplications of the generic elliptic fiber. Since the involution $-1$ is always there, the only possible groups of complex multiplications of elliptic curves have order $2,4$ or $6$. Hence $k \in \{1,2,3,4,6\}$.

Now $\t_B = t_\x \circ \a_B$, where $\x := \t_B(\s)$. Since $\a_B^m = 1$, we have that $(\t_B)^m$ is simply given by translation by a section. Since it is of finite order, it must be given by translation by a torsion section. Hence, the order $d$ of $(\t_B)^m$ divides the order of the torsion subgroup $\MWt$.
\end{proof}

\begin{lemma} \label{l:fixed}
If $k>1$, then $(\tau_B) ^ {\bar m}$ has a fixed curve intersecting each elliptic fiber.
\end{lemma}

\begin{proof}
As we have seen, if $k > 1$, $(\a_B)^{\bar m}$ is given by complex multiplication on the generic elliptic fiber, which has a fixed curve (the zero section) intersecting each elliptic fiber. $(\tau_B) ^ {\bar m}$ simply composes $(\a_B)^{\bar m}$ with translation by a section; therefore, $(\tau_B) ^ {\bar m}$ also has a fixed curve intersecting each elliptic fiber.
\end{proof}

Let us now define the group homomorphism $\cp_i : MW \to MW$ by
\begin{equation}\label{e:cp}
\cp_i: \x \mapsto \cp_i(\x) := \a_B^{i-1} (\x) \boxplus \a_B^{i-2} (\x) \boxplus \ldots \boxplus \x,
\end{equation}
where $\x \in MW$ and $\boxplus$ denotes addition of sections in the Mordell-Weil group.

Define also the map $\PH_i : MW \to \IQ$ by
\begin{equation}\label{e:phi}
\PH_i: \x \mapsto \PH_i(\x) := \langle \cp_i(\x), \cp_i(\x) \rangle,
\end{equation}
where $\x \in MW$ and $\langle \cdot, \cdot \rangle$ is the height pairing introduced in \eqref{e:hp}.

\begin{lemma} \label{l:data}
The following data are equivalent:
\begin{enumerate}
\item $\tau_B \in \Aut(B)$, of finite order $n$;
\item Pairs $(\alpha_B,\xi)$, with $\alpha_B \in \Aut_\s(B)$ of order $m$ 
dividing $n$
and
$\Phi_m(\xi)=0$
(or equivalently, $\cp_m(\xi) \in \MWt$,
a torsion section of order $d=n/m$).
\end{enumerate}
\end{lemma}

\begin{proof}

Let $\t_B \in \Aut(B)$ be of finite order $n$, with linearization $\a_B := \l(\t_B) \in \Aut_\s(B)$ of order $m$. Let $\x = \t_B(\s)$; then $\t_B = t_\x \circ \a_B$. By definition, $\a_B(\s) = \s$. We need to show that $\PH_m(\x)=0$.

We know that 
\begin{align}
1 =& \t_B^n \notag\\
=& (t_\x \circ \a_B)^n  \notag\\
=&  (t_{\cp_m (\x)} \circ \a_B^m)^d \notag\\
=& (t_{\cp_m (\x)})^d,
\end{align}
where we used the fact that $\a_B^m = 1$. This implies that $\cp_m(\x)$ is a torsion section of order $d=n/m$. From proposition \ref{p:torsion}, we know that a section $\h$ is torsion if and only if $\langle \h, \h \rangle =0$. Hence $\PH_m(\x) = 0$. 

Conversely, given an automorphism $\a_B \in \Aut_\s(B)$ of order $m$, and a section 
$\x$
with $\Phi_m(\xi)=0$, 
we form the automorphism $\t_B = t_\x \circ \a_B$. Then $\t_B$ is of order $n$, with $n = m d$ where $d$ is the order of the torsion section $\cp_m(\x)$.
\end{proof}

It follows that in order to classify all automorphisms $\t_B \in \Aut(B)$ of finite order, we must construct all pairs $(\a_B, \x)$ as above. 

\subsubsection{Automorphisms of the first kind of finite order}\label{s:tors}
An automorphism $\t_B \in \Aut(B)$ is of the first kind if and only if $\bar m =1$, hence 
the order of $\a_B$ is $m = k$.
From the lemmas above, either $\t_B$ has a fixed curve intersecting each elliptic fiber, or $\a_B$ is the identity and $\t_B = t_\h$, where $t_\h$ means translation by a torsion section $\h$ of order $n=d$.

The complete list of rational elliptic surfaces with finite subgroups of the automorphism groups generated by translations by torsion sections is easy to produce. 
All possible configurations of singular fibers of rational elliptic surfaces
are classified in \cite{Miranda:1990, Persson:1990}, and for each configuration, 
the torsion group $\MWt$ of the Mordell-Weil group of sections is tabulated. 
This was also described in terms of lattices in \cite{Oguiso:1991}. Hence, the complete list of finite subgroups of the automorphism groups generated by translations by torsion sections can be harvested naturally from the list of torsion groups of rational elliptic surfaces.



\subsubsection{Automorphisms of the second kind of finite order}\label{s:p1}

Automorphisms of the second kind correspond to the cases where $\bar m > 1$. Then $\a_B$ is of order $m = \bar m k > 1$. From 
the lemmas above, either $(\t_B)^{\bar m}$ 
has a fixed curve intersecting each elliptic fiber, or $k=1$ and $m = \bar m > 1$, \emph{i.e.} the order of $\a_B$ is equal to the order of the induced automorphism $\t_{\IP^1}$ on the $\IP^1$ base.

\subsection{Restrictions}\label{s:restrictions}

In this work, we are interested in automorphisms of rational elliptic surfaces that can be lifted to free automorphisms on smooth fiber products of two rational elliptic surfaces. Accordingly, we impose two restrictions: the fiber product must be smooth, and the automorphism on the fiber product must be free. Since we focus on automorphisms of rational elliptic surfaces, we would like to understand what these restrictions on the fiber product impose on the constituent automorphisms of rational elliptic surfaces.

Take two rational elliptic surfaces $B$ and $B'$, and form the fiber product $\tX = B \times_{\IP^1} B'$. 

Suppose that $\t_B \in \Aut(B)$ has integers $(n,m,\bar m,k,d)$ defined in the previous section. Similarly, suppose 
that $\t_{B'} \in \Aut(B')$ has integers $(n',m',\bar m', k', d')$. One can form the automorphism 
$\t_\tX = \t_B \times_{\IP^1} \t_{B'}$ on the fiber product $\tX$ if and only if they have consistent induced action 
on the $\IP^1$: ${\rho}(\tau_{B}) = {\rho}'(\tau_{B'})$
and in particular, $\bar m = \bar m'$. 
In fact, any automorphism $\t_\tX \in \Aut(\tX)$ that preserves the elliptic fibration of $\tX$ can be constructed in this way.

We want the fiber product $\tX$ to be smooth, and the group $\langle \t_\tX \rangle$ generated by $\t_\tX$ to act freely on $\tX$. Recall from the introduction that $\tX$ is smooth if and only if $S'' =0$, that is, the singular fibers of $B$ and $B'$ are not paired in the fiber product. 
Furthermore, it is clear that $\t_\tX$ will be free if and only if the fibers of $B$ and $B'$ intersecting the fixed locus of $\t_B$ and $\t_{B'}$ are not paired in the fiber product $\tX = B \times_{\IP^1} B'$.

But what do these restrictions imply on the rational elliptic surfaces $B$ and $B'$ and the automorphisms $\t_B$ and $\t_{B'}$?

\begin{lemma}
If $\langle \t_\tX \rangle$ is free on $\tX$, then $n=n'$ and $k=k'=1$.
\end{lemma}

\begin{proof}
Suppose that $n<n'$. Then, $(\t_B)^{n}$ is the identity on $B$; hence, $(\t_\tX)^n$ is not free on the fiber product $B \times_{\IP^1} B'$. Similarly for $n' < n$; therefore if $\langle \t_\tX \rangle$ is free, then $n=n'$.

Suppose now that $k > 1$. Lemma \ref{l:fixed} tells us that $(\t_B)^{\bar m}$ has a fixed curve intersecting each elliptic fiber of $B$. Since the fixed locus of $(\t_{B'})^{\bar m}$ cannot be empty, some of the fibers of $B$ and $B'$ intersecting the fixed locus of $(\t_B)^{\bar m}$ and $(\t_{B'})^{\bar m}$ will necessarily be paired in the fiber product. Hence, $(\t_\tX)^{\bar m}$ will not be free. Similarly for $k' > 1$. Therefore, if $\langle \t_\tX \rangle$ is free, then $k=k'=1$.
\end{proof}

Henceforth we may set $n=n'$ and $k = k' = 1$, which implies that $d=d'$ and $m = \bar m = \bar m' = m'$, that is, the orders of $\a_B$ and $\a_{B'}$ are the same as the order of the induced automorphism $\t_{\IP^1}$ of the $\IP^1$ base. We are left with only two independent integers, $(n, m)$, with $n = d m$; both $\t_B$ and $\t_{B'}$ have order $n$, and $\a_B$, $\a_{B'}$ and $\t_{\IP^1}$ all have order $m$.

\subsubsection{$m=1$}
Consider first the case with $m=1$ (automorphisms of the first kind). In this case, both $\t_B$ and $\t_{B'}$ are given by translation by torsion sections $\h$ and $\h'$ of order $n=d$, and there is no induced action on the $\IP^1$. These automorphisms act freely on the smooth elliptic fibers of $B$ and $B'$ respectively. Hence, the only fibers of $B$ (resp. $B'$) intersecting the fixed locus of $\t_B$ (resp. $\t_{B'}$) must be singular. By choosing the isomorphism between the $\IP^1$ bases of $B$ and $B'$ to form the fiber product, in this case we can always make sure that the singular fibers of $B$ and $B'$ are not paired in the fiber product, thus ensuring simultaneously that the fiber product $\tX$ is smooth and the automorphism group $\langle \t_{\tX} \rangle$ is free. 

\subsubsection{$m>1$}
Suppose now that $m>1$ (automorphisms of the second kind). Without loss of generality, we can choose homogeneous coordinates $[u:v]$ on $\IP^1$ such that $\t_{\IP^1}$ is given by
\begin{equation}\label{e:p1map}
\t_{\IP^1}: [u:v] \mapsto [u: \o v],
\end{equation}
where $\o$ is an $m$'th root of unity. 

\begin{lemma}\label{l:mg1}
Suppose that $m>1$. $\langle \t_\tX \rangle$ is free if and only if $\langle \t_B \rangle $ is free on $f_\infty$ and $\langle \t_{B'} \rangle$ on $f'_0$, or vice-versa. Moreover, since we assume that $\tX$ is smooth, $f_0$ and $f'_0$ are not both singular, and the same holds for $f_\infty$ and $f'_\infty$.
\end{lemma}

\begin{proof}
Using the parameterization above, it is clear that $0,\infty \in \IP^1$ are the only fixed points of $\langle \t_{\IP^1} \rangle$. 
The fixed locus of $\t_B$ (resp. $\t_{B'}$) must then be contained 
in the fibers $f_0 := \b^{-1}(0)$ (resp. $f_0':=(\b')^{-1}(0)$) and $f_\infty:= \b^{-1}(\infty)$ (resp. $f_\infty':=(\b')^{-1}(\infty)$). Moreover, we know that the fixed loci of $\t_B$ and $\t_{B'}$ cannot be empty; hence, $\t_B$ cannot be free on $f_0$ and $f_\infty$ at the same time; the same holds for $\t_{B'}$. Therefore, the automorphism group generated by $\t_\tX= \t_B \times_{\IP^1} \t_{B'}$ is free on $\tX$ if and only if $\langle \t_B \rangle $ is free on $f_\infty$ and $\langle \t_{B'} \rangle$ on $f'_0$, or vice-versa.

The second claim is that if $\tX$ is smooth, then $f_0$ and $f'_0$ are not both singular, and the same holds for $f_\infty$ and $f'_\infty$. The fibers of $\tX$ above $0$ and $\infty$ are given by $f_0 \times f'_0$ and $f_\infty \times f'_\infty$ respectively. Hence, the claim follows directly from the fact that $\tX$ is smooth if and only if $S''=0$.
\end{proof}

\subsubsection{$m>1$, $d>1$}
We can extract a stronger restriction for the case $m>1$, $d>1$, due to the following lemma.

\begin{lemma}\label{l:d1}
For $m>1$ and $d>1$, if $\tau_B$ is free on $f_\infty$, then $f_0$ must be
singular.
\end{lemma}

\begin{proof}
Suppose that $m>1$, $d>1$, and $\t_B$ is free on $f_\infty$. 
Since $d>1$, $(\t_B)^m$ is given by translation by a non-zero torsion section of order $d$. 
Since $\t_B$ is free on $f_\infty$, it cannot be free on $f_0$, so
$(\t_B)^m$ is also not free on $f_0$. However, translation by a non-zero torsion section acts freely on a smooth elliptic fiber.
Therefore $f_0$ must be singular.
\end{proof}

Hence, we get the following corollary.

\begin{corollary}
Suppose that $m>1$, $d>1$, $\langle \t_\tX \rangle$ is free, and $\tX$ is smooth. 
Then, up to interchanging $0$ and $\infty$ if necessary:
\begin{itemize} 
\item $\langle \t_B \rangle $ is free on $f_\infty$ and $\langle \t_{B'} \rangle$ on $f'_0$;
\item $f_0$ and $f'_\infty$ are singular;
\item $f_\infty$ and $f'_0$ are smooth.
\end{itemize}
\end{corollary}

\begin{proof}
This follows immediately from lemmas 3.6 and 3.7.
\end{proof}

We now have quite a few restrictions on the rational elliptic surfaces $B$ and $B'$, and their automorphisms $\t_B$ and $\t_{B'}$, which we formalize in the following subsection.

\subsubsection{Suitable pairs}

From the previous analysis, we introduce the concept of a suitable pair $(B,\t_B)$.

\begin{definition}\label{d:suitable}
Let $\t_B \in \Aut(B)$ with integers $(n,m,\bar m, d, k)$. We say that the pair $(B,\t_B)$ is  \emph{suitable} if $k=1$, which implies that $m = \bar m$, and one of the following conditions holds: 
\begin{enumerate}
\item $m=1$ (first kind);
\item $m >1$ (second kind), $d=1$, and the group generated by $\t_B$ acts freely on $f_\infty$;
\item $m > 1$ (second kind), $d> 1$, $f_\infty$ is smooth, and the group generated by $\t_B$ acts freely on $f_\infty$. 
(In this case lemma \ref{l:d1} implies that $f_0$ is singular.)
\end{enumerate}
\end{definition}

According to the discussion above, any free action on a smooth fiber
product comes from a pair of suitable pairs $(B, \tau_B)$ and $(B', \tau_{B'})$, and conversely any
suitable pair $(B, \tau_B)$ can be lifted to a free automorphism on a
smooth fiber product by combining it with another appropriate suitable pair $(B',\tau_{B'})$.

Notice that when $m=1$ (condition 1), 
suitable pairs are already classified by torsion groups of rational elliptic surfaces, as explained in section \ref{s:tors}. Our main goal is then to classify suitable pairs $(B,\t_B)$ satisfying 
conditions 2 or 3.

\begin{proposition}\label{p:finf}
Let $(B,\t_B)$ be a suitable pair with $m>1$, $d=1$ (condition 2). Then $f_\infty$ is either smooth or singular of type $I_{n r}$, for some integer $r > 0$.
\end{proposition}

\begin{proof}
By definition, $f_\infty$ is either a smooth elliptic curve or a singular fiber of one of Kodaira's type. By looking at Kodaira's list of singular fibers \cite{Kodaira}, it is easy to see that the only type of singular fiber that admits a free cyclic automorphism of order $n$ is $I_{nr}$ for some 
integer $r>0$.
\end{proof}

Hence, when $m>1$, $d=1$ (condition 2), $f_\infty$ must be either smooth or of type $I_{nr}$, while when $m>1$, $d>1$ (condition 3), $f_\infty$ must be smooth.

We can also extract a requirement on the linearizations $\a_B$.

\begin{proposition}\label{p:abinf}
Let $(B,\t_B)$ be a suitable pair with $m>1$ (conditions 2 and 3). Then the linearization $\a_B$ must fix the neutral component of the fiber $f_\infty$ pointwise.
\end{proposition}

\begin{proof}
Suppose first that 
$f_\infty$ is smooth. Since $(B,\t_B)$ is suitable, $\t_B$ acts freely on $f_\infty$. We know that $\a_B$ must fix the zero section; hence it cannot be free on $f_\infty$. Therefore, if $\t_B = t_\x \circ \a_B$ is free on $f_\infty$, then $\a_B$ must fix $f_\infty$ pointwise.

Suppose now that 
$f_\infty$ 
is of type $I_{nr}$. Let $\TH_0, \ldots, \TH_{nr-1}$ be the components of the singular fibers $f_\infty$. Since $\a_B$ fixes the zero section, it must send the neutral component $\TH_0$ to itself. Now $\a_B$ must also sends each other components $\TH_i$, $i=1, \ldots, nr-1$ to itself, otherwise $\t_B  = t_\x \circ \a_B$ would have fixed points. Hence $\a_B$ must send each component to itself. Then, for $\t_B$ to be free, on each component $\a_B$ must either fix the component pointwise or fix only the intersection points of the component with its neighbor components. But since $\a_B$ fixes the zero section, which is not an intersection point of the neutral component $\TH_0$ with its neighbors, we obtain that $\a_B$ must fix the neutral component of $f_\infty$ pointwise.
\end{proof}

We are now in a position to classify suitable pairs $(B,\t_B)$ satisfying conditions 2 and 3. Using the insight of lemma \ref{l:data}, we first study the linearizations $\a_B$, and then scrutinize the subgroup of sections $\ker(\PH_m) \subseteq MW$.

\section{The linearizations $\a_B$}\label{s:linear}

Let us study first the linearizations $\a_B$. 

As usual, let $\a_B \in \Aut_\s(B)$ be of order $m$ and $\t_{\IP^1} \in \Aut(\IP^1)$ be of order $\bar m$. Recall from the previous section that we want to classify all $B$ admitting an $\a_B$ such that
\begin{enumerate}
\item $m = \bar m$;
\item $f_\infty$ of $B$ is either smooth or of type $I_{mr}$;
\item $\a_B$ fixes the neutral component of $f_\infty$ pointwise.
\end{enumerate}

\begin{definition}\label{d:sigma}
We will call a pair $(B,\a_B)$ satisfying these three conditions a \emph{suitable $\s$-pair} (not to be confused with suitable pairs $(B,\t_B)$).
\end{definition}
Since $\a_B$ fixes the zero section it must act on the Weierstrass model of $B$. We could use this point of view to classify all such automorphisms. However, we prefer to use a different viewpoint. 

Suppose that $(B,\a_B)$ is a suitable $\s$-pair. Let $\cb = B / \langle
\a_B \rangle$ be the quotient surface. Let $\widehat{\cb}$ be the relatively minimal (or Kodaira) model of $\cb$.
Recall that this is obtained from a resolution ${\widehat{\cb}}^+$ of
$\cb$ by blowing down all $(-1)$-curves in fibers. Although
$\widehat{\cb}$ may no longer map to $\cb$, it does map to $\IP^1$, in
fact it is an elliptic fibration over $\IP^1$. Since it is birational to a
quotient of a rational surface, it is itself rational, so it is a smooth
rational elliptic surface.

Hence, any suitable $\s$-pair $(B,\a_B)$ with $\a_B$ of order $m$ fits in a commutative diagram
\begin{equation}\label{e:quot}
\xymatrix{
\widehat{B} \ar[r] \ar[d] & \widehat{\cb} \ar@{.>}[d] \ar@/^1.5pc/[dd]\\
B \ar[d]^{\b} \ar[r] & \cb \ar[d] \\
\IP^1 \ar[r]^g & \IP^1
}
\end{equation}
Here $g: \IP^1 \to \IP^1$ is the $m$-th power map totally ramified at
$0,\infty \in \IP^1$, and $\widehat{B}$ is the minimal resolution of the
fiber product ${\widehat{B}}^{-}$.

Conversely, we can start with any rational elliptic surface
$\widehat{\cb}$ and recover $(B,\a_B)$. First we pull back via the map
$g$ to get the fiber product ${\widehat{B}}^{-}$. The $\IP^1$ which is the
source of $g$ has an automorphism $\alpha_{\IP^1}$ of order $m$ which is a deck transformation with respect to $g$, that is $g \circ \alpha_{\IP^1} = g$ (and generates the group of such automorphisms); this together with the
identity on $\widehat{\cb}$ induce a natural automorphism
$\a_{{\widehat{B}}^{-}}$ on the fiber product ${\widehat{B}}^{-}$. We
resolve ${\widehat{B}}^{-}$ to get $\widehat{B}$; and then blow down the
$(-1)$ curves to obtain a $\s$-pair $(B,\alpha_B)$. In general, uniqueness
of the Kodaira model guarantees that starting with $(B,\alpha_B)$, going
to the Kodaira model $\widehat{\cb}$ of its quotient, and then going
through the above steps brings us back to the original $(B,\alpha_B)$. We
need to understand which smooth rational elliptic surfaces $\widehat{\cb}$
are suitable in the sense that they yield a suitable $\s$-pair
$(B,\alpha_B)$.

We illustrate the issue with two examples. First, consider a $\widehat{\cb}$
whose fiber over $\infty$ is of type $I_M$. The inverse image in the fiber
product ${\widehat{B}}^{-}$ of the original $I_M$ fiber still looks like an $I_M$
fiber, but the surface ${\widehat{B}}^{-}$ in which it sits is singular: it has
an $A_{m-1}$ singularity at each of the $M$ singular points of the $I_M$ fiber.
We resolve these singularities to obtain the smooth $\widehat{B}$, containing an
$I_{mM}$ fiber over $\infty$. This elliptic surface $\widehat{B}$ is minimal,
so is isomorphic to $B$. The natural automorphism of the fiber product
${\widehat{B}}^{-}$ induces an $\alpha_B$ on $B$. This $\alpha_B$ acts trivially on
the $M$ original components, but non-trivially (via various roots of unity)
on the $(m-1)M$ new components arising from the resolution of the
$A_{m-1}$ singularities. In particular, we see that we have a suitable
sigma-pair $(B,\alpha_B)$.

On the other hand, consider the case that the fiber of $\widehat{\cb}$ over
$\infty$ is of type $I_0^*$, and $m=2$. The fiber product ${\widehat{B}}^{-}$ is now
not normal. When we normalize, the central component of the $I_0^*$ fiber is
replaced with an elliptic curve. The resulting elliptic surface $\widehat{B}$ is
not minimal; $B$ is obtained from it by blowing down the four rational
components of the $\infty$ fiber, and has a smooth ($I_0$) fiber there.  In
this case $(B,\alpha_B)$ is not a suitable sigma-pair even though its fiber
over $\infty$ is smooth: $\alpha_B$ acts non-trivially on this smooth elliptic
curve which is the neutral (and only) component.

The general pattern is obtained in Lemma 4.6 below, and depends on the
notion of deficiency which we now define.

\subsection{Deficiencies}

The first step consists in understanding the relationship between the singular fibers of $B$ and $\widehat{\cb}$. If the map $g$ is etale at $p \in \IP^1$ --- which in this case is true for any $p \in \IP^1$ not equal to $0$ or $\infty$ --- then the fiber of $\widehat{\cb}$ (or $\cb$) above $g(p)$ is identical to the fiber of $B$ above $p$. The crucial analysis is then for the fibers $f_0$ and $f_\infty$ of $B$, where the map $g$ is totally ramified of order $m$.

This question was addressed by Miranda and Persson in section 7 of \cite{Miranda:1986} using the Weierstrass models of $\widehat{\cb}$ and $B$. Table 7.1 of \cite{Miranda:1986} gives the relation between the types of the fiber of $\widehat{\cb}$ over $g(p)$ and the fiber of $B$ over $p$ for the points $p$ where $g$ is ramified of order $m$.\footnote{In their notation, $B$ is $Y$ and $\widehat{\cb}$ is $X$.} We reproduce this table in table \ref{t:ram} for convenience.

\begin{table}[!hbt]
\begin{center}
\begin{footnotesize}
\begin{tabular}{|c|lll|}
\hline
Fiber of $\widehat{\cb}$ over $g(p)$ & Fiber of $B$ over $p$ && \\
\hline
$I_0$ & $I_0$ && \\[0.3em]
$I_M$ & $I_{m M}$ && \\[0.3em]
$I_M^*$ & $I_{m M}$ if $m$ even; & $I_{m M}^*$ if $m$ odd &\\[0.3em]
$II$&$I_0$ if $m =0$ mod $6$; & $II$ if $m=1$ mod $6$; & $IV$ if $m=2$ mod $6$\\
&$I_0^*$ if $m =3$ mod $6$; & $IV^*$ if $m=4$ mod $6$; & $II^*$ if $m=5$ mod $6$\\[0.3em]
$III$ & $I_0$ if $m =0$ mod $4$; & $III$ if $m=1$ mod $4$; & \\
&$I_0^*$ if $m =2$ mod $4$; & $III^*$ if $m=3$ mod $4$; & \\[0.3em]
$IV$&$I_0$ if $m =0$ mod $3$; & $IV$ if $m=1$ mod $3$; & $IV^*$ if $m=2$ mod $3$\\[0.3em]
$IV^*$&$I_0$ if $m =0$ mod $3$; & $IV^*$ if $m=1$ mod $3$; & $IV$ if $m=2$ mod $3$\\[0.3em]
$III^*$ & $I_0$ if $m =0$ mod $4$; & $III^*$ if $m=1$ mod $4$; & \\
&$I_0^*$ if $m=2$ mod $4$; & $III$ if $m=3$ mod $4$; & \\[0.3em]
$II^*$&$I_0$ if $m =0$ mod $6$; & $II^*$ if $m=1$ mod $6$; & $IV^*$ if $m=2$ mod $6$\\
&$I_0^*$ if $m =3$ mod $6$; & $IV$ if $m=4$ mod $6$; & $II$ if $m=5$ mod $6$\\
\hline
\end{tabular}
\end{footnotesize}\vspace{0.5em}
\caption{Relation between the types of singular fibers of $B$ and $\widehat{\cb}$ over the ramified points $p$ of order $m$ ($p=0,\infty$) \cite{Miranda:1986}.}\label{t:ram}
\end{center}
\end{table}

The next step consists in understanding the allowed configurations of singular fibers for $\widehat{\cb}$.

We can define an integer, which we call the \emph{deficiency}, associated to singular fibers of $\widehat{\cb}$, and depending on the order $m$ of $g$. Basically, the deficiency formalizes the change in Euler characteristic of the singular fibers between $B$ and $\widehat{\cb}$.

\begin{definition}
Let $\cs$ be a fiber of $\widehat{\cb}$ over a point $p \in \IP^1$. We define the \emph{deficiency} $D(\cs)$ of $\cs$ by
\begin{equation}\label{e:def}
D(\cs) = \c(\cs) - \sum_{S \in W} \c(S),
\end{equation}
where $W$ is the set of fibers over the preimages $g^{-1}(p)$ of $p$. 
\end{definition}

\begin{proposition}
$D(\cs) = 0$ if $\cs$ is smooth.
\end{proposition}

\begin{proof}
This is obvious from \eqref{e:def}, table \ref{t:ram} and the fact that the fiber of $\widehat{\cb}$ above $g(p)$ is identical to the fiber of $B$ above $p$ at the points $p$ where $g$ is etale.
\end{proof}

\begin{lemma}\label{l:sum0}
Let $\car$ be the set of singular fibers of $\widehat{\cb}$. Then
$$
\sum_{\cs \in \car} D(\cs) = 0.
$$
\end{lemma}

\begin{proof}
By definition,
$$
\sum_{\cs \in \car} D(\cs) = \left( \sum_{\cs \in \car} \c(\cs)\right) - \left( \sum_{\cs \in \car} \sum_{S \in S_\cs} \c(S)\right).
$$
But $\sum_{\cs \in \car} \c(\cs) = \c(\widehat{\cb}) = 12$, and since any singular fiber of $B$ is pullbacked from a singular fiber of $\widehat{\cb}$,
\begin{align}
\sum_{\cs \in \car} \sum_{S \in S_\cs} \c(S) &= \sum_{S \in R} \c(S)\notag\\
 &=12,\notag
\end{align}
where $R$ is the set of singular fibers of $B$. Hence, $\sum_{\cs \in \car} D(\cs) = 0$.
\end{proof}

Using table \ref{t:ram} relating singular fibers of $\widehat{\cb}$ and $B$ over the ramified points $0,\infty \in \IP^1$, we can find explicit formulae for the deficiencies of the singular fibers of $\widehat{\cb}$.

\begin{lemma}\label{l:def}
Let $\cs$ be a singular fiber of $\widehat{\cb}$ over a point $p \in \IP^1$. Then the deficiency of $\cs$ is given as in table \ref{t:def}.
\begin{table}[!htb]
\begin{center}
\begin{footnotesize}
\begin{tabular}{|lll|}
\hline
Position of $p$ & Type of $\cs$ & Deficiency $D(\cs)$\\
\hline
$p\neq 0,\infty$ & Any type & $D(\cs) = (1-m)\c (\cs)$\\[0.3em]
$p=0$ & $I_M$ & $D(\cs) = (1-m) \c(\cs)$\\
 or $p=\infty$&$I_M^*$ & $D(\cs) = (1-m) \c(\cs) + 6 (m - \d)$ with $\d = m ~\rm{mod}~2$\\
&$II,III,IV$ & $D(\cs) = (1-\d )\c(\cs)$ with $\d = m~\rm{mod}~\left( 12/\c(\cs) \right)$\\
&$II^*,III^*,IV^*$ & $D(\cs) = \epsilon (12-\c(\cs))$ with $\epsilon = (m-1) ~\rm{mod}~\left( 12/(12 - \c(\cs)) \right)$\\
\hline
\end{tabular}
\end{footnotesize}
\vspace{0.5em}	
\caption{Deficiencies of singular fibers $\cs$ of $\cb$.}\label{t:def}
\end{center}
\end{table}

\end{lemma}

\begin{proof}
This is a simple calculation from the definition of the deficiency \eqref{e:def} and table \ref{t:ram}. Note that the last line cannot be expressed neatly in terms of $\delta$, due
to the peculiarity of clock arithmetic: for $\delta \neq 0$, we have
$\epsilon=\delta-1$, but this fails for $\delta=0$.
\end{proof}

\subsection{Suitable rational elliptic surfaces}

Using lemma \ref{l:sum0} and table \ref{t:def} we are in position to generate a complete list of suitable quotient surfaces $\widehat{\cb}$, taking into account the three requirements in the definition of a suitable $\s$-pair $(B,\a_B)$ in the beginning of this section.

\begin{lemma}\label{l:cond}
Let $\a_B \in \Aut_\s(B)$ be of order $m$. If $(B,\a_B)$ is a suitable $\s$-pair, then the fiber $f_\infty$ of $\widehat{\cb}$ is either smooth or of type $I_r$ for some $r\geq 1$, and
$$
D(f_0) = (m-1)\left(12 - \c(f_0)\right),
$$
where $f_0$ is the fiber above $0$ of $\widehat{\cb}$. In this case we say that $\widehat{\cb}$ is \emph{suitable of order $m$}.
\end{lemma}

\begin{proof}
If $(B,\a_B)$ is a suitable $\s$-pair, then $f_\infty$ of $B$ is either smooth or of type $I_{mr}$, and $\a_B$ fixes the neutral component of $f_\infty$ pointwise.

Suppose that $f_\infty$ is smooth. Since $\a_B$ fixes $f_\infty$ pointwise, it is clear that the fiber $f_\infty$ of $\widehat{\cb}$ will also be smooth.\footnote{In general, the fixed locus of any automorphism $\alpha$ of a
manifold is a submanifold \cite{Cartan}. The differential $d{\alpha}$ acting on the
tangent space at a fixed point has the tangent space to the fixed locus as
its $+1$ eigenspace.  The quotient is smooth if and only if there is a
single normal direction, \emph{i.e.} $\text{codim}=1$. 
In our case, we have a curve of fixed points on a smooth surface, so the
quotient is smooth there, hence the fiber (in the resolved $\widehat{\cb}$) is
the naive quotient, which is a smooth elliptic curve.}

If $f_\infty$ of $B$ is of type $I_{m r}$, then table \ref{t:ram} says that $f_\infty$ of $\widehat{\cb}$ is either of type $I_{r}$ or of type $I^*_{r}$ and $m$ is even. But the latter case occurs only if $\a_B$ does not fix the neutral component of $f_\infty$ of $B$ pointwise. Hence $f_\infty$ of $\widehat{\cb}$ must be of type $I_r$.

Now, let $\car$ be the set of singular fibers of $\widehat{\cb}$. We know that $\sum_{\cs \in \car} D(\cs) =0$. Using table \ref{t:def}, and the fact that $f_\infty$ of $\widehat{\cb}$ is either smooth or of type $I_r$, we get
$$
0=\sum_{\cs \in \car} D(\cs) = D(f_0) + (1-m) \left (\c(f_\infty)+\sum_{\cs \in \car \setminus \{f_0,f_\infty\}} \c(\cs) \right). 
$$
But for a rational elliptic surface
$$
\c(f_\infty)+ \sum_{\cs \in \car \setminus \{f_0,f_\infty\}}\c(\cs)  = 12 - \c(f_0),
$$  
and the lemma is proved.
\end{proof}

We can now produce a list of all suitable quotient surfaces $\widehat{\cb}$.

\begin{lemma}\label{l:suit}
A rational elliptic surface $\widehat{\cb}$ is suitable of order $m$ if and only if $f_\infty$ is either smooth or of type $I_r$ and $f_0$ is of one of the types listed in table \ref{t:f0}. All other $m$'s are not allowed.
\begin{table}[!hbt]
\begin{center}
\begin{tabular}{|c|c|}
\hline
$m$ & Type of $f_0$\\
\hline
$6$ & $II^*$\\
$5$ & $II^*$\\
$4$ & $II^*, III^*$\\
$3$ & $II^*, III^*, IV^*$\\
$2$ & $II^*, III^*, IV^*,I_4^*, I_3^*, I_2^*, I_1^*, I_0^*$\\
\hline
\end{tabular}\vspace{0.5em}
\caption{List of allowed $f_0$ for $\widehat{\cb}$}\label{t:f0}
\end{center}
\end{table}
\end{lemma}

\begin{proof}
This follows from lemma \ref{l:cond} and table \ref{t:def} after a straightforward calculation.
\end{proof}

This is enough to produce a complete list of suitable quotient surfaces $\widehat{\cb}$. We scan through Persson's list of configurations of singular fibers \cite{Persson:1990}, and keep all configurations with $f_0$ as in lemma \ref{l:suit} and $f_\infty$ either smooth or of type $I_k$. Each of these configurations is suitable for one or more $m$'s.

\subsection{List of suitable $\s$-pairs $(B,\a_B)$}

To produce the list of suitable $\s$-pairs $(B,\a_B)$, we pullback each rational elliptic surface $\widehat{\cb}$ suitable of order $m$ that we found in the previous section, via the map $g: \IP^1 \to \IP^1$ of order $m$ as in \eqref{e:quot}. We use table \ref{t:ram} to extract the types of the fibers $f_0$ and $f_\infty$ of $B$. This generates all suitable $\s$-pairs $(B,\a_B)$.

\subsubsection{A `simplification'}\label{s:simple}

The choice of the $m$-sheeted covering map $g$ (as in \eqref{e:quot}) is equivalent to the choice of its two branch points $0,\infty \in \IP^1$. Now the point $0$ is rigid,
since by lemma \ref{l:suit} it must correspond to a singular fiber. But we are
free to move the point $\infty$ over the Zariski open subset of $\IP^1$
parametrizing fibers of type $I_r$ (including the smooth ones, $r=0$). So
any $\widehat{\cb}$ and $m$ as in lemma \ref{l:suit} produce a $1$-parameter family of suitable
$\s$-pairs $(B,\alpha_B)$. At most finitely many of these involve a
singular fiber $f_\infty$. It follows that any suitable 
$\s$-pair $(B,\alpha_B)$ coming from a $\widehat{\cb}$ with singular fiber $f_\infty$ is the
smooth specialization of a $1$-parameter family of suitable $\s$-pairs
$(B,\alpha_B)$ coming from quotient surfaces $\widehat{\cb}$ with smooth fibers
$f_\infty$. 

From now on we will therefore focus on surfaces $B$ and $\widehat{\cb}$
with $f_\infty$ smooth.

\subsubsection{The list}

Keeping this in mind, we get the following result.

\begin{proposition}
The list of suitable $\widehat{\cb}$ and their associated suitable $\s$-pairs $(B,\a_B)$ with smooth $f_\infty$ is given in table \ref{t:pairs}. We indicate the Mordell-Weil lattices of the rational elliptic surfaces $B$ and $\widehat{\cb}$: for each Mordell-Weil lattice the allowed configurations of singular fibers can be found using Oguiso and Shioda's and Persson's list \cite{Miranda:1990,Oguiso:1991,Persson:1990}. We also write down the root lattice $T$ associated to the singular fibers of $B$, as in \cite{Oguiso:1991}. In the table $I_0$ means a smooth fiber.
\end{proposition}

Note that while writing down this table we found a small mistake in Oguiso and Shioda's list \cite{Oguiso:1991} for case 30 in our list, the quotient surface of which corresponds to case 33 in their list. There, it is stated that $MW = A_1^* \oplus \langle 1/6 \rangle$ for the quotient surface $\widehat{\cb}$, but it should be $MW =\frac{1}{6} \begin{pmatrix} 2 & 1\\1 & 2\end{pmatrix}$.

\begin{table}[!hbt]
\begin{center}
\begin{footnotesize}
\begin{tabular}{|c|ccc||cccc|c|}
\hline
$m$&\multicolumn{3}{|c||}{$\widehat{\cb}$} & \multicolumn{4}{c|}{$B$}&$num$ \\
\hline
 & $f_0$ & $\MWl$ & $\MWt$ & $f_0$&$\MWl$&$\MWt$ &$T$& \\
\hline
$6$ & $II^*$ & $0$ & $0$ & $I_0$ & $E_8$ & $0$ &$0$&  $1$\\
\hline
$5$ & $II^*$ & $0 $ & $0$ & $II$ & $E_8$ & $0$ &$0$& $2$\\
\hline
$4$ & $II^*$ & $0 $ & $0$ & $IV$  & $E_6^*$ & $0$ &$A_2$& $3$ \\
& $III^* $ & $A_1^* $ & $0 $ & $ I_0$ & $E_8$ & $0$ &$0$& $4$ \\
&& $0$ & $\IZ_2$ &  & $D_4^*$ & $\IZ_2$ & ${A_1}^{\oplus 4}$&$5$ \\
\hline
$3$ & $II^*$ & $0 $ & $0$ & $I_0^*$  & $D_4^*$ & $0$ &$D_4$& $6$ \\
& $III^* $ & $A_1^* $ & $0 $ & $ III$ & $E_7^*$ & $0$&$A_1$ & $7$ \\
&& $0$ & $\IZ_2$ &  & $D_4^*$ & $\IZ_2$ & ${A_1}^{\oplus 4}$&$8$ \\
&$IV^*$ & $A_2^*$ & $0$&$I_0$&$E_8$&$0$&$0$&$9$\\
&& $\langle 1/6 \rangle$ & $0$&&$D_4^* \oplus A_1^*$&$0$&${A_1}^{\oplus 3}$&$10$\\
&& $0$ & $\IZ_3$&&$A_2^*$&$\IZ_3$&${A_2}^{\oplus 3}$&$11$\\
\hline
$2$ & $II^*$ & $0 $ & $0$ & $IV^*$  & $A_2^*$ & $0$ &$E_6$& $12$ \\
& $III^* $ & $A_1^* $ & $0 $ & $ I_0^*$ & $D_4^*$ & $0$ &$D_4$& $13$ \\
&& $0$ & $\IZ_2$ &  & ${A_1^*}^{\oplus 2}$ & $\IZ_2$ &$D_4 \oplus {A_1}^{\oplus 2}$& $14$ \\
&$IV^*$ & $A_2^*$ & $0$&$IV$&$E_6^*$&$0$&$A_2$&$15$\\
&& $\langle 1/6 \rangle$ & $0$&&$\frac{1}{6}\begin{pmatrix}2 & 1 &0&-1\\1&5&3&1\\0&3&6&3\\-1&1&3&5 \end{pmatrix}$&$0$&$A_2 \oplus {A_1}^{\oplus 2}$&$16$\\
&& $0$ & $\IZ_3$&&$A_2^*$&$\IZ_3$&${A_2}^{\oplus 3}$&$17$\\
&$I_4^*$ & $0$ & $\IZ_2 $ & $I_8$ & $A_1^* $ & $\IZ_2 $ &$A_7$& $18$\\
&$I_3^*$ & $\langle 1/4 \rangle$ & $0$ & $I_6$ & $A_2^* \oplus A_1^*$ & $0$ &$A_5$& $19$\\
&$I_2^*$ & ${A_1^*}^{\oplus 2}$ & $0$ & $I_4$ & $D_5^*$ & $0$ &$A_3$& $20$\\
&&$A_1^*$ & $\IZ_2$ &  & $A_3^*$ & $\IZ_2$ &$A_3 \oplus {A_1}^{\oplus 2}$& $21$ \\
&&$0$ & $(\IZ_2)^2$ &  & $\langle 1/4 \rangle$ & $(\IZ_2)^2$ & $A_3 \oplus {A_1}^{\oplus 4}$&$22$\\
&$I_1^*$&$A_3^*$ & $0$ & $I_2$ & $E_7^*$ & $0$&$A_1$&$23$\\
&&$A_1^* \oplus \langle 1/4 \rangle$& $0$ & & $D_4^* \oplus A_1^*$ &$0$ &${A_1}^{\oplus 3}$& $24$\\
&&$\langle 1/12 \rangle$ & $0$ & & $A_2^* \oplus \langle 1/6 \rangle $ & $0$ & ${A_2}^{\oplus 2} \oplus A_1$&$25$\\
&&$\langle 1/4 \rangle$ & $\IZ_2$ & & ${A_1^*}^{\oplus 3} $ & $\IZ_2$ & ${A_1}^{\oplus 5}$&$26$\\
&&$0$&$\IZ_4$ & & $A_1^*$ & $\IZ_4$ &${A_3}^{\oplus 2} \oplus A_1$& $27$\\
&$I_0^*$ & $D_4^*$ & $0$ & $I_0$ & $E_8$ & $0$ & $0$&$28$\\
&&$ {A_1^*}^{\oplus 3}$ & $0$ & & $D_6^*$ & $0$ & ${A_1}^{\oplus 2}$&$29$\\
&&$\frac{1}{6} \begin{pmatrix} 2 & 1\\1 & 2\end{pmatrix}$ & $0$ & & ${A_2^*}^{\oplus 2}$& $0$ &${A_2}^{\oplus 2}$& $30$\\
&&${A_1^*}^{\oplus 2}$ & $\IZ_2$ & & $D_4^*$ & $\IZ_2$ &${A_1}^{\oplus 4}$& $31$\\
&&$A_1^*$ & $(\IZ_2)^2$ & & ${A_1^*}^{\oplus 2}$ & $(\IZ_2)^2$ &${A_1}^{\oplus 6}$& $32$\\
&&$\langle 1/4 \rangle$ & $\IZ_2 $ & & $ {A_1^*}^{\oplus 2}$ & $\IZ_2$ &${A_3}^{\oplus 2}$& $33$\\
&&$0$&$(\IZ_2)^2$ && $0$ & $(\IZ_2)^2$ &${D_4}^{\oplus 2}$& $34$\\
\hline
\end{tabular}
\end{footnotesize}\vspace{0.5em}
\caption{List of suitable $\widehat{\cb}$ and the associated suitable $\s$-pairs $(B,\a_B)$ with smooth $f_\infty$.}\label{t:pairs}
\end{center}
\end{table}

\section{The sections $\x$}\label{s:sections}

In the previous section we produced a list of suitable $\s$-pairs $(B,\a_B)$. We now want to complete the pairs $(B,\a_B)$ to suitable pairs $(B,\t_B)$ as defined in definition \ref{d:suitable}. More precisely, we want to complete the automorphisms $\a_B$ to automorphisms $\t_B = t_\xi \circ \a_B$ of the second kind satisfying conditions 2 or 3 in definition \ref{d:suitable}, where $\x \in MW$ is a section.

Throughout this section we let $m, n = m d$ denote the orders of $\alpha_B,
\tau_B$ respectively. We recall from lemma 3.4 that $\xi$ must satisfy
$\Phi_m(\xi)=0$. Equivalently, $\cp_m(\xi)$ is a torsion section; its order
is precisely the integer $d=n/m$. This is also the order of the point
where $\cp_m(\xi)$ meets the fiber $f_\infty$, which by the `simplification' of
subsection 4.3.1 can and will be taken to be smooth.

We will determine the group $\ker(\Phi_m)$ in subsection 5.1. In 5.2 we
will determine its subgroup $\ker(\cp_m) = {\ker(\Phi_m)}_{d=1}$; the
possibility that $d>1$ arises when these two groups are not equal. The
remaining condition for suitability is that the action on $f_{\infty}$
should be free; this is discussed in 5.3. An overview of the entire
analysis is given in 5.4.

\subsection{Calculation of $\ker(\PH_m)$}



\begin{lemma}\label{l:all}
$\ker(\PH_m) = \left[(MW^{\a_B})^\perp\text{ in }MW\right]$.\footnote{$(MW^{\a_B})^\perp$ in $MW$ denotes the orthogonal complement of $MW^{\a_B}$ in $MW$ with respect to the height pairing.}
\end{lemma}

\begin{proof}
Consider the map $MW \to MW \otimes \IC$. Under this map, $\alpha_B$ becomes an automorphism of order $m$ of a vector space, and so its
roots are $m$-th roots of $1$.  But $\PH_m$ kills all roots other than $1$, so
$(\PH_m)/m$ is projection onto the invariants, and its kernel is the subspace perpendicular to the space of $\a_B$-invariants.  But both $\ker(\PH_m)$ and $\left[(MW^{\a_B})^\perp\text{ in }MW\right]$ are the
inverse images of their images under $MW \to MW \otimes \IC$; therefore $\ker(\PH_m) = \left[(MW^{\a_B})^\perp\text{ in }MW\right]$.
\end{proof}

We now know how to find $\ker(\PH_m)$; we only have to calculate, for each case, the orthogonal complement of $MW^{\a_B}$ in $MW$, and we know that $MW^{\a_B} = m \widehat{MW}$,\footnote{By $m A$ we mean that the intersection matrix of $A$ is multiplied by $m$.} where $\widehat{MW}$ is the Mordell-Weil group of the quotient surface $\widehat{B}$ and $m$ is the order of $\a_B$ \cite{Shioda:1990}.

This can be done for each case, using the following result.

\begin{lemma}\label{l:perps}
For each of the following cases, the embedding is unique up to isomorphisms, and the orthogonal complement is as follows. We write in brackets the case number of table \ref{t:pairs} that these embeddings correspond to. In all the other cases of table \ref{t:pairs}, $\MWl(\widehat{\cb}) = 0$.
\begin{align}
&\left[ A_1^\perp\text{ in }E_8\right]= E_7~~~~(4), &&\left[ \langle 3/2 \rangle^\perp\text{ in }E_7^*\right] = E_6^* ~~~~(7), \notag\\
&\left[ A_2^\perp\text{ in }E_8\right] = E_6 ~~~~(9), &&\left[ (A_1^*)^\perp\text{ in }D_4^* \oplus A_1^*\right] = D_4^* ~~~~(10), \notag\\
&\left[ (\langle 1 \rangle)^\perp\text{ in }D_4^*\right] = \langle 1 \rangle^{\oplus 3}~~~~ (13), &&\left[ (2 A_2^*)^\perp\text{ in }E_6^*\right] = D_4 ~~~~(15), \notag\\
&\left[ (\langle 1/3 \rangle)^\perp\text{ in }\frac{1}{6}\begin{pmatrix}2 & 1 &0&-1\\1&5&3&1\\0&3&6&3\\-1&1&3&5 \end{pmatrix} \right] = \langle 1 \rangle^{\oplus 3}~~~~ (16), \notag\\
&\left[ (A_1^*)^\perp\text{ in }A_2^* \oplus A_1^*\right] = A_2^* ~~~~(19),&
&\left[ ({\langle 1 \rangle}^{\oplus 2})^\perp\text{ in }D_5^*\right] = \langle 1 \rangle^{\oplus 3} ~~~~(20),\notag\\
&\left[ (\langle 1 \rangle \oplus \IZ_2)^\perp\text{ in }A_3^* \oplus \IZ_2 \right] = \langle 1 \rangle^{\oplus 2} \oplus \IZ_2 ~~~~(21), 
&&\left[ (2 A_3^*)^\perp\text{ in }E_7^*\right] = D_4~~~~ (23),\notag\\
&\left[ (\langle 1 \rangle \oplus A_1^*)^\perp\text{ in }D_4^* \oplus A_1^*\right] = \langle 1 \rangle^{\oplus 3}~~~~ (24), &&\left[ (\langle 1/6 \rangle)^\perp\text{ in }A_2^* \oplus \langle 1/6 \rangle \right] = A_2^*~~~~ (25),\notag\\
&\left[ (A_1^* \oplus \IZ_2)^\perp\text{ in }{A_1^*}^{\oplus 3} \oplus \IZ_2 \right] = {A_1^*}^{\oplus 2} \oplus \IZ_2 ~~~~(26), &&\left[ (D_4)^\perp\text{ in }E_8 \right] = D_4 ~~~~(28), \notag\\
&\left[ ({\langle 1 \rangle}^{\oplus 3})^\perp\text{ in }D_6^* \right]=\langle 1 \rangle^{\oplus 3} ~~~~(29),
&&\left[ (A_2^*)^\perp\text{ in }{A_2^*}^{\oplus 2} \right] = A_2^*~~~~ (30),\notag\\
&\left[ ({\langle 1 \rangle}^{\oplus 2} \oplus \IZ_2)^\perp\text{ in }D_4^*\oplus \IZ_2 \right]  ={\langle 1 \rangle}^{\oplus 2} \oplus \IZ_2~~~~ (31), \notag\\
&\left[ ({\langle 1 \rangle} \oplus (\IZ_2)^2)^\perp\text{ in }{A_1^*}^{\oplus 2}  \oplus (\IZ_2)^2 \right]={\langle 1 \rangle} \oplus (\IZ_2)^2 ~~~~(32),\notag\\
&\left[ (A_1^* \oplus \IZ_2)^\perp\text{ in }{A_1^*}^{\oplus 2} \oplus \IZ_2 \right] = A_1^* \oplus \IZ_2~~~~ (33). \notag
\end{align}
\end{lemma}

\begin{proof}
For most of the cases this is just a simple exercise in linear algebra. Here we only prove a few cases; the other cases are either straighforward or can be proved similarly.

\emph{Case 32.} We want $\left[ ({\langle 1 \rangle} \oplus (\IZ_2)^2)^\perp\text{ in }{A_1^*}^{\oplus 2}  \oplus (\IZ_2)^2 \right]$. There are only 4 elements of ${A_1^*}^{\oplus 2}$ of squared length 1: the sums of one generator from each $A_1^*$.  In ${A_1^*}^{\oplus 2}  \oplus (\IZ_2)^2$, we can use the same 4 elements plus the torsion elements. The embedding of $\langle 1 \rangle \oplus (\IZ_2)^2$ is determined by the choice of one of these 4 elements of ${A_1^*}^{\oplus 2}$ plus the torsion elements. The orthogonal complement is the sum of the entire torsion $(\IZ_2)^2$ plus the orthogonal complement in ${A_1^*}^{\oplus 2}$, which is the complementary copy of $\langle 1 \rangle$.

\emph{Case 31.} Here we want $\left[ ({\langle 1 \rangle}^{\oplus 2} \oplus \IZ_2)^\perp\text{ in }D_4^*\oplus \IZ_2 \right]$. Let us use the explicit description of $D_4^*$ as the square lattice $\IZ^4$ plus the single non-integral element $\frac{1}{2}(1,1,1,1)$. The 24 elements of length 1 are plus/minus the 4 unit vectors in $\IZ^4$, and $\frac{1}{2}(\pm 1,\pm 1,\pm 1,\pm 1)$. Triality interchanges all of these, so we can take the image of the generator of the first $\langle 1 \rangle$ to be, say, $(1,0,0,0)$. The image of the generator of the second $\langle 1 \rangle$ must then be $(0,1,0,0)$, up to permutation. The orthogonal complement is then ${(0,0,a,b)}=\langle 1 \rangle^2$, and $\left[ ({\langle 1 \rangle}^{\oplus 2} \oplus \IZ_2)^\perp\text{ in }D_4^*\oplus \IZ_2 \right]={\langle 1 \rangle}^{\oplus 2} \oplus \IZ_2$.

\emph{Case 15.} We want $\left[ (2 A_2^*)^\perp\text{ in }E_6^*\right]$. That is, we need two vectors in $E_6^*$ with length $4/3$ and intersection $-2/3$. The minimal points in $E_6^*$ have length $4/3$, and there are $54$ of them: $30$ of the form $\pm \frac{1}{3}(0,-2,-2,1,1,1,1,0)$, and another $24$ of the form $\pm \frac{1}{6} (\pm 3,5,-1,-1,-1,-1,-1,\mp 3)$ \cite{Conway}. We find three possible embeddings, given by sending the two generators of $A_2^*$ to one of the three following possibilities, up to permutations:
\begin{gather}
\{\frac{1}{3}(0,-2,-2,1,1,1,1,0), \frac{1}{3}(0,1,1,-2,-2,1,1,0)\}\notag\\
\{\frac{1}{3}(0,-2,-2,1,1,1,1,0), \frac{1}{6}(3,5,-1,-1,-1,-1,-1,-3)\}\notag\\
\{\frac{1}{6}(-3,-1,5,-1,-1,-1,-1,3), \frac{1}{6}(3,5,-1,-1,-1,-1,-1,-3)\}\notag
\end{gather}
These three embeddings are presumably equivalent, but we did not check
this. Instead, we perform the calculations separately in the three
cases.  In each case, we calculate the intersection numbers of these vectors with the basis of $E_6^*$ given in \cite{Conway}, p.127, eq.(125). From this we work out the orthogonal complement for each embedding in this basis of $E_6^*$, and we get in each case the intersection matrix
$$
\begin{pmatrix}
2&-1&0&0\\
-1&2&-1&-2\\
0&-1&2&0\\
0&-1&0&2
\end{pmatrix}
$$
for the orthogonal complement, which corresponds to $D_4$.

\emph{Case 23.} We want to compute $\left[ (2 A_3^*)^\perp\text{ in }E_7^*\right]$. The intersection matrix of $2 A_3^*$ is
$$\frac{1}{2}
\begin{pmatrix}
3&2&1\\
2&4&2\\
1&2&3
\end{pmatrix}.
$$
So we are looking for three elements $a_1$, $a_2$ and $a_3$ of $E_7^*$ with these intersection numbers. The minimal norm in $E_7^*$ is $3/2$, and there are 56 minimal vectors of the form $\pm \frac{1}{4} (-3,-3,1,1,1,1,1,1)$ \cite{Conway}. So without loss of generality we can take $a_1 = \frac{1}{4} (-3,-3,1,1,1,1,1,1)$. Now $a_3$ must also be a vector of minimal length, and it must have intersection number $\frac{1}{2}$ with $a_1$. There are two possibilities, up to permutations: $a_3 = \frac{1}{4}(-3,1,-3,1,1,1,1,1)$ or $a_3'=\frac{1}{4}(-1,-1,-1,-1,-1,-1,3,3)$. Now $a_2$ must have norm $2$. We note that vectors in $E_7$ have integer norm, while vectors in the nontrivial coset $E_7^*-E_7$ have half-integer (but not integer) norm. Hence, $a_2$ must be in $E_7$; that is, it must be a minimal vector of $E_7$. It is a permutation of either $(-1,1,0,0,0,0,0,0)$ or $\frac{1}{2}(1,1,1,1,-1,-1,-1,-1)$. Imposing the required intersection numbers with $a_3$ and $a_1$, we get two possibilities for $a_2$, regardless of whether we consider $a_3$ or $a_3'$: $a_2 = (-1,0,0,0,0,0,0,1)$ or $a_2' = \frac{1}{2}(-1,-1,-1,-1,1,1,1,1)$, up to permutations. Hence in total we have four choices of embeddings for $2 A_3^*$ in $E_7^*$. Then, we compute the intersection numbers of our $a$'s with the generators of $E_7^*$ given in \cite{Conway}, p. 125, eq.(115). It is simple linear algebra to find the orthogonal complement of our four embeddings in $E_7^*$ in this basis, and it is straightforward to show that all four orthogonal complements have intersection matrix  
$$
\begin{pmatrix}
2&-1&0&0\\
-1&2&-1&-2\\
0&-1&2&0\\
0&-1&0&2
\end{pmatrix},
$$
which again corresponds to $D_4$.

\emph{Case 16.} In this case we want $\left[ \langle 1/3 \rangle^\perp\text{ in }L^* \right]$, where $L = \left[ A_1^\perp\text{ in } A_5 \right]$. The intersection matrix for $L^*$ was computed in \cite{Oguiso:1991}; this is the matrix in the statement of the lemma. Now we want to embed a vector of length $\frac{1}{3}$ in $L^*$. Note that the quotient $L^*/L$ is $(\IZ_2)^2 \times \IZ_3$. We write down representatives for the 12 nontrivial cosets of $L$ in $L^*$. The square of each element of $L^*$ is a rational number which modulo $2$ depends only on its image in $L^*/L$. The possible values are $0, \frac{1}{3}, \frac{5}{6}, 1, \frac{3}{2}, \frac{4}{3} \text{ mod } 2$. It follows that the shortest possible length
squared is $\frac{1}{3}$, and any minimal vector is, up to equivalence, of the form $\pm \frac{1}{6} (1,1,1,1,-2,-2)$. Taking an explicit basis for $L^*$ we compute the intersection number of the generators with our minimal vector, and find the orthogonal complement in $L^*$ in this basis. We find that the intersection matrix is simply the $3 \times 3$ identity matrix, hence the orthogonal complement is the cubic lattice $\IZ^3$, or in the above notation $\langle 1 \rangle^{\oplus 3}$.

\end{proof}

Combining lemmas \ref{l:all} and \ref{l:perps} we have an explicit lattice description of $\ker(\PH_m)$ for each case in table \ref{t:pairs}.

\subsection{Cases where $d>1$ is allowed}

The next step consists in finding the cases where $d>1$ is allowed. Given a section $\x \in \ker(\PH_m)$, we know that the automorphism $\t_B = t_\x \circ \a_B$ is of order $n$, $n = d m$ for an integer $ d \geq 1$. However, the integer $d$ depends on the section $\x$. For some sections $\xi \in \ker(\PH_m)$, $d$ may equal $1$, while for others it may be greater than $1$. We need to find out when this happens.



We showed in lemma \ref{l:k} that $d$ divides the order of $\MWt$. Hence, if $\MWt=0$, then $d=1$. We also showed in lemma \ref{l:d1} that if $f_0$ is smooth, then $d=1$. Thus, the only cases in table \ref{t:pairs} that may have $d > 1$ are cases 8, 14, 17, 18, 21, 22, 26 and 27.

For the cases where $d > 1$ is possible, the group $\ker(\Phi_m)$ will split into the disjoint union of a subgroup $\ker(\Phi_m)_{d=1}$ of sections yielding $d=1$, and its complement $\ker(\Phi_m)_{d>1}$ of sections with $d>1$. 

\begin{lemma}\label{l:d}
The only cases where $d>1$ is allowed are listed in table \ref{t:d}. We write down the case number of table \ref{t:pairs}, the allowed integer $d$, the group $\ker(\Phi_m)$, and its subgroup $\ker(\Phi_m)_{d=1}$. For each case with $d=2$, $\ker(\Phi_m)_{d=2}$ is the non-trivial coset of $\ker(\Phi_m)_{d=1}$ in $\ker(\Phi_m)$, while when $d=3$, $\ker(\Phi_m)_{d=3}$ is the union of the two non-trivial cosets of $\ker(\Phi_m)_{d=1}$ in $\ker(\Phi_m)$.
\begin{table}[!htb]
\begin{center}
\begin{tabular}{|c|ccc|}
\hline
Case &$d$&$\ker(\Phi_m)$&$\ker(\Phi_m)_{d=1}$\\
\hline
8 & $2$ & $D_4^* \oplus \IZ_2$&$D_4^*$\\
14 & $2$ & ${A_1^*}^{\oplus 2} \oplus \IZ_2$ & $\langle 1 \rangle^{\oplus 2} \oplus \IZ_2$\\
17 & $3$ & $A_2^* \oplus \IZ_3$&$A_2^*$\\
22 & $2$ & $\langle 1/4 \rangle \oplus (\IZ_2)^2$&$\langle 1 \rangle \oplus (\IZ_2)^2$\\
26 & $2$ & ${A_1^*}^{\oplus 2} \oplus \IZ_2$ & $\langle 1 \rangle^{\oplus 2} \oplus \IZ_2$\\
27 & $2$ & $A_1^* \oplus \IZ_4$ & $A_1^* \oplus \IZ_2$\\
\hline
\end{tabular}
\vspace{0.5em}
\caption{List of cases where $d>1$ is allowed.}\label{t:d}
\end{center}
\end{table}
\end{lemma}

\begin{proof}
We will prove this lemma case by case. But let us first state the general philosophy of the proof. For the 8 possible cases listed before the lemma, given a section $\x \in \ker(\PH_m)$, we want to know whether $\cp_m(\x)$ is the zero section or a non-zero torsion section. Our strategy goes as follows. We first find the action of $\a_B$ on the components of the singular fibers. Then we use proposition \ref{p:torsion} and the fact that the torsion sections are $\a_B$-invariant to find the components that the torsion sections intersect. Finally, using this result and the fact that the height pairing of any section with a torsion section is zero, we are able to decide whether $d>1$ is allowed, and find the decomposition of $\ker(\Phi_m)$.\\

\noindent\emph{Case 8.} $m=3$. $\ker(\PH_3) = MW = D_4^* \oplus \IZ_2$. $f_0$ is of type $III$, and $T = {A_1}^{\oplus 4}$. Let $R_i$, $i=0,1$ denote the components of $f_0$, and $S_{j,k}$, with $j=1,2,3$ and $k=0,1$ the components of the three other singular fibers with associated root lattice $A_1$. The zero components are the neutral components as will always be the case in the following.

$\a_B$ acts on the components as follows. On $f_0$, it sends each component to itself: $\a_B : R_i \mapsto R_i$ for $i=0,1$. $\a_B$ permutes the other singular fibers, hence if we denote them appropriately we obtain that $\a_B : S_{j,k} \mapsto S_{j+1,k}$ for $j=1, 2, 3$ and $k=0,1$, where $S_{4,k} \equiv S_{1,k}$ is understood.

Let $\h$ be the non-zero torsion section of order $2$. We know from proposition \ref{p:torsion} that $\sum_{s_i} {\rm contr}_{s_i} (\h,\h) = 2$, where $s_i$ denote the reducible fibers of $B$. For a singular fiber $s$ with associated root lattice $A_1$, we have that ${\rm contr}_s (\h, \h) = 0$ if $\h$ intersects the neutral component, and ${\rm contr}_s (\h, \h) = 1/2$ otherwise. Hence, we directly obtain that $\h$ must intersect $R_1$ and $S_{j,1}$ for $j=1,2,3$.

Take a section $\x \in \ker(\PH_3)$. Suppose that it intersects $R_0$. Then, from the action of $\a_B$ on $f_0$ we get that $\cp_3(\x)$ also intersects $R_0$, hence $\cp_3(\x) = \s$ and $d=1$. Now suppose that $\x$ intersects $R_1$; we get that $\cp_3(\x)$ also intersects $R_1$, that is $\cp_3(\x) = \h$ and $d=2$. Therefore, any section $\x \in \ker(\PH_3)$ intersecting $R_0$ satisfies $\cp_3(\x) = \s$, while any section $\x \in \ker(\PH_3)$ intersecting $R_1$ satisfies $\cp_3(\x) = \h$.

Now take any section $\x \in \ker(\PH_3)$. If it intersects $R_0$, then $\x \boxplus \h$\footnote{As usual $\boxplus$ denote addition in the Mordell-Weil group.} intersects $R_1$; and vice-versa if $\x$ intersects $R_1$. Hence, we obtain that $\ker(\PH_3)$ breaks into the disjoint union of its index 2 sublattice $\ker(\Phi_3)_{d=1} = D_4^*$, and the non-trivial coset $\ker(\Phi_3)_{d=2}$.\\

\noindent\emph{Case 14.} $m=2$. $\ker(\PH_2) = MW = {A_1^*}^{\oplus 2} \oplus \IZ_2$. $f_0$ is of type $I_0^*$, and $T=D_4 \oplus {A_1}^{\oplus 2}$. We denote the five irreducible components of $f_0$ by $R_{0,0}$, $R_{1,0}$, $R_{0,1}$, $R_{1,1}$, and $R'$; $R_{0,0}$ is the neutral component, and $R'$ is the component of $I_0^*$ with multiplicity $2$. We write $S_{j,k}$, with $j=1,2$ and $k=0,1$ for the components of the two other singular fibers with root lattices $A_1$.

$\a_B$ acts on the two $A_1$ fibers by $\a_B : S_{j,k} \mapsto S_{j+1,k}$, with $S_{3,k} \equiv S_{1,k}$. Using the fact that the fiber $f_0$ of the quotient surface $\widehat{\cb}$ is of type $III^*$, we obtain the action of $\a_B$ on the components of $f_0$:
$$
\a_B:~~\{ R_{0,0}, R_{1,0}, R_{0,1}, R_{1,1}, R'\} \mapsto \{ R_{0,0}, R_{0,1}, R_{1,0}, R_{1,1}, R' \}.
$$

The contributions ${\rm contr}(\h, \h)$ associated to fibers with root lattices $A_1$ are either $0$ or $1/2$. For a fiber of type $I_0^*$ (with root lattice $D_4$), the contributions can be either $0$ for the neutral component or $1$ for the four components with subscripts --- recall that no sections intersect the component $R'$ with multiplicity $2$. Hence, from proposition \ref{p:torsion}, and the fact that the non-zero torsion section $\h$ of order $2$ must be $\a_B$-invariant, we get that $\h$ intersects the non-neutral components $R_{1,1}$, $S_{1,1}$ and $S_{2,1}$.

From the action of $\a_B$, we get that any section $\x \in \ker(\PH_2)$ intersecting either $R_{0,0}$ or $R_{1,1}$ must satisfy $\cp_2(\x) = \s$ and $d=1$, while any section $\x \in \ker(\PH_2)$ intersecting either $R_{0,1}$ or $R_{1,0}$ must satisfy $\cp_2(\x) = \h$ and $d=2$. However, we need more information to determine the intersection numbers of the sections $\x \in \ker(\PH_2)$.

Now let $\x \in \ker(\PH_2)$ be a section corresponding to a minimal point of the Mordell-Weil lattice ${A_1^*}^{\oplus 2}$, which has length $1/2$. Since for this configuration of singular fibers $\sum_{s_i} {\rm contr}_{s_i} (x) \leq 2$, we know that $\x$ must be disjoint from the zero section. Hence
$$
\langle \x, \x \rangle = 1/2 = 2 - \sum_{s_i} {\rm contr}_{s_i} (\x).
$$
Using this constraint, the fact that $\langle \x, \h \rangle = 0$ and the intersection numbers of $\h$, we can list the possible intersection numbers of $\x$. We obtain 4 possibilities: $(R_{1,0}, S_{1,1}, S_{2,0})$, $(R_{1,0}, S_{1,0}, S_{2,1})$, $(R_{0,1}, S_{1,0}, S_{2,1})$ and $(R_{0,1}, S_{1,1}, S_{2,0})$. In all of these cases, $\x$ intersects either $R_{0,1}$ or $R_{1,0}$. Hence, the sections $\x$ corresponding to the minimal points of the Mordell-Weil lattice must satisfy $\cp_2(\x) = \h$.

Therefore, we get that $\ker(\PH_2)$ splits into the disjoint union of the 45 degree rotated sublattice $\ker(\PH_2)_{d=1} = {\langle 1 \rangle}^{\oplus 2} \oplus \IZ_2$ and its non-trivial coset $\ker(\PH_2)_{d=2}$.\\

\noindent\emph{Case 17.} $m=2$. $\ker(\PH_2) = MW = A_2^* \oplus \IZ_3$. $f_0$ is of type $IV$, and $T = {A_2}^{\oplus 3}$. As usual, let $R_i$, $i=0,1,2$ be the components of $f_0$, and $S_{j,k}$, $j=1,2$ and $k=0,1,2$ be the components of the two other singular fibers.

The contributions associated to these three reducible fibers with root lattices $A_2$ are either $0$ for the neutral component, or $2/3$ for the two other irreducible components. 

$\a_B$ sends each component of $f_0$ to itself, that is $\a_B : R_i \mapsto R_i$ for $i=0,1,2$. For the two other fibers, as usual $\a_B: S_{j,k} \mapsto S_{j+1,k}$ for $j=1,2$ and $k=0,1,2$, with $S_{3,k} \equiv S_{1,k}$.

There are two non-zero torsion sections of order $3$; denote them by $\h$ and $\h \boxplus \h$. Using proposition \ref{p:torsion}, we find that $\h$ intersects $R_1$, $S_{1,1}$ and $S_{2,1}$, and accordingly $\h \boxplus \h$ intersects $R_2$, $S_{1,2}$ and $S_{2,2}$. 

Since $\a_B$ sends each component of $f_0$ to itself, we get that any section $\x \in \ker(\PH_2)$ intersecting the neutral component $R_0$ must satisfy $\cp_2(\x)=\s$, any section $\x \in \ker(\PH_2)$ intersecting $R_1$ must satisfy $\cp_2(\x)=\h \boxplus \h$, while any section $\x \in \ker(\PH_2)$ intersecting $R_2$ must satisfy $\cp_2(\x)=\h$. In the two last cases, we obtain that $d=3$, since $\h$ and $\h \boxplus \h$ are of order $3$.

Now take any section $\x \in \ker(\PH_2)$. Suppose that $\x$ intersects $R_l$ for some $l$. Then, from the intersection form of the torsion sections, we get that $\x \boxplus \h$ will intersect $R_{l+1}$, and $\x \boxplus \h \boxplus \h$ will intersect $R_{l+2}$, where $R_3 \equiv R_0$. 

Hence, we get that $\ker(\PH_2)$ breaks into the disjoint union of the index 3 sublattice $\ker(\PH_2)_{d=1} = A_2^*$, and $\ker(\PH_2)_{d=3}$ which is the union of its two non-trivial cosets.\\

\noindent\emph{Case 18.} $m=2$. $\ker(\PH_2) = MW = A_1^* \oplus \IZ_2$. $f_0$ is of type $I_8$, and $T = A_7$. Let $R_i$, $i=0, \ldots, 7$ denote the components of $f_0$.

The contributions for $f_0$ are $0$ for the neutral component, $7/8$ for $R_1$ and $R_7$, $3/2$ for $R_2$ and $R_6$, $15/8$ for $R_3$ and $R_5$, and $2$ for $R_4$. Hence, from proposition \ref{p:torsion} the non-zero torsion section $\h$ of order $2$ must intersect $R_4$.

The action of $\a_B$ on the components of $f_0$ is given by $\a_B: R_i \mapsto R_{8-i}$, for $i = 0, \ldots, 7$, with $R_8 \equiv R_0$. Hence, we directly get that for any section $\x \ker(\PH_2)$, $\cp_2(\x)$ must intersect $R_0$, that is $\cp_2(\x) = \s$ and $d=1$. Therefore, it turns out that in this case $d > 1$ is not possible.\\

\noindent\emph{Case 21.} $m=2$. $\ker(\PH_2) = \langle 1 \rangle^{\oplus 2} \oplus \IZ_2$.  $f_0$ is of type $I_4$, and $T = A_3 \oplus {A_1}^{\oplus 2}$. Denote by $R_i$, $i=0,1,2,3$ the components of $f_0$, and by $S_{j,k}$, $j=1,2$ and $k=0,1$ the components of the two other singular fibers with root lattices $A_1$.

The contributions for $f_0$ are either $0$ for the neutral component, $3/4$ for $R_1$ and $R_3$, and $1$ for $R_2$. For the two other reducible fibers, it is either $0$ for the neutral components, or $1/2$ otherwise.

From proposition \ref{p:torsion}, the non-zero torsion section $\h$ of order $2$ must intersect $R_2$, $S_{1,1}$ and $S_{2,1}$. Since the action of $\a_B$ on $f_0$ is given by $\a_B: R_i \mapsto R_{4-i}$ for $i=0,1,2,3$, with $R_4 \equiv R_0$, we get that for any section $\x \in \ker(\PH_2)$, $\cp_2(\x)$ must intersect $R_0$, that is $\cp_2(\x)= \s$. Therefore, in this case as well it turns out that $d>1$ is not possible. \\

\noindent\emph{Case 22.} $m=2$. $\ker(\PH_2) = MW = \langle 1/4 \rangle \oplus (\IZ_2)^2$. $f_0$ is of type $I_4$, and $T=A_3 \oplus {A_1}^{\oplus 4}$. Denote by $R_i$, $i=0,1,2,3$ the components of $f_0$, and by $S_{j,k}$, $j=1,2,3,4$ and $k=0,1$ the components of the four other singular fibers with root lattices $A_1$.

The contributions for $f_0$ are either $0$ for the neutral component, $3/4$ for $R_1$ and $R_3$, and $1$ for $R_2$. For the four other reducible fibers, it is either $0$ for the neutral components, or $1/2$ otherwise.

The action of $\a_B$ on $f_0$ is given by $\a_B: R_i \mapsto R_{4-i}$ for $i=0,1,2,3$, with $R_4 \equiv R_0$. On the four other fibers, it is given by
$$
\a_B:~~\{ S_{1,k},S_{2,k},S_{3,k},S_{4,k} \} \mapsto \{ S_{2,k},S_{1,k},S_{4,k},S_{3,k} \},
$$
for $k=0,1$.

There are three non-zero torsion sections of order $2$, which we denote by $\h_1$, $\h_2$ and $\h_1 \boxplus \h_2$. From proposition \ref{p:torsion} and the fact that they are all $\a_B$-invariant, we get that $\h_1$ intersects the non-neutral components $R_2$, $S_{1,1}$ and $S_{2,1}$, $\h_2$ intersects $R_2$, $S_{3,1}$ and $S_{4,1}$, and $\h_1 \boxplus \h_2$ intersects $S_{1,1}$, $S_{2,1}$, $S_{3,1}$ and $S_{4,1}$.

Take any section $\x \in \ker(\PH_2)$. From the action of $\a_B$ on $f_0$, we know that $\cp_2(\x)$ must intersect $R_0$. Hence, the only possibilities are that $\cp_2(\x)=\s$ or $\cp_2(\x) = \h_1 \boxplus \h_2$. To distinguish between these two cases we must look at the intersection numbers of the sections $\x$ with the other singular fibers.

Let $\x \in \ker(\PH_2)$ correspond to a minimal point in the Mordell-Weil lattice $\langle 1/4 \rangle$, which has length $1/4$. Since for this configuration of singular fibers $\sum_{s_i} {\rm contr}_{s_i} (\x) \leq 3$, we know that $\x$ must be disjoint from $\s$. Hence
$$
\langle \x, \x \rangle = 1/4 = 2 - \sum_{s_i} {\rm contr}_{s_i} (\x).
$$
We can use this constraint, the fact that $\langle \x, \h_1 \rangle = \langle \x, \h_2 \rangle = \langle \x, \h_1 \boxplus \h_2 \rangle = 0$ and the intersection numbers of the torsion sections above to list the possible intersection numbers of $\x$. We get the 8 different combinations (listing only the non-neutral components that $\x$ intersects) $(R_{1 \text{ or } 3}, S_{1,1 \text{ or } 2,1}, S_{3,1 \text{ or } 4,1})$. Using the action of $\a_B$ above, for each of these cases it is easy to see that the section $\cp_2(\x)$ will intersect the components $R_0, S_{1,1}, S_{2,1}, S_{3,1}$ and $S_{4,1}$. That is, for any section $\x$ corresponding to a minimal point of the Mordell-Weil lattice, we have that $\cp_2(\x) = \h_1 \boxplus \h_2$ and $d=2$.

Therefore, we obtain that $\ker(\PH_2)$ splits into the disjoint union of the sublattice $\ker(\PH_2)_{d=1} = \langle 1 \rangle \oplus (\IZ_2)^2$ and its non-trivial coset $\ker(\PH_2)_{d=2}$.\\

\noindent\emph{Case 26.} $m=2$. $\ker(\PH_2) = {A_1^*}^{\oplus 2} \oplus \IZ_2$. $f_0$ is of type $I_2$, and $T = {A_1}^{\oplus 5}$. Let $R_i$, $i=0,1$ be the components of $f_0$, and $S_{j,k}$, $j=1,2,3,4$, $k=0,1$ be the components of the four other singular fibers with root lattices $A_1$.

The contributions for all the singular fibers are either $0$ for the neutral component, or $1/2$ otherwise. $\a_B$ sends each component of $f_0$ to itself, and permutes two by two the other components. That is, if we name the components appropriately,
$$
\a_B:~~\{ R_i, S_{1,k}, S_{2,k}, S_{3,k}, S_{4,k} \} \mapsto \{ R_i, S_{2,k}, S_{1,k}, S_{4,k}, S_{3,k} \},
$$
for $i=0,1$ and $k=0,1$. 

There is one non-zero torsion section of order $2$, which we denote by $\h$. From proposition \ref{p:torsion} and the fact that $\h$ is $\a_B$-invariant, we get that it must intersect the components $R_0, S_{1,1}, S_{2,1}, S_{3,1}$ and $S_{4,1}$.

Take any section $\x \in \ker(\PH_2)$. From the action of $\a_B$ on $f_0$ we know that $\cp_2(\x)$ must intersect $R_0$. Hence, both possibilities, $\cp_2(\x) = \s$ or $\cp_2(\x) = \h$, may occur.

Suppose now that $\x \in MW$ corresponds to a minimal point in the Mordell-Weil lattice $\MWl = {A_1^*}^{\oplus 3}$, hence it has length $1/2$. Since for this configuration $\sum_{s_i} {\rm contr}_{s_i}(\x, \x) \leq 5/2$, it follows that $\x$ is disjoint from the zero section. Hence we get
$$
\langle \x, \x \rangle = 1/2 = 2 - \sum_{s_i} {\rm contr}_{s_i}(\x, \x).
$$
Using also the fact that $\langle \x, \h \rangle = 0$ and the intersection numbers of $\h$, we can list the possible intersection numbers of the minimal $\x$. We obtain 6 possibilities: listing only the non-neutral components that $\x$ intersects, $(R_1, S_{1,1}, S_{2,1})$, $(R_1, S_{3,1}, S_{4,1})$, $(R_1, S_{1,1}, S_{3,1})$, $(R_1, S_{1,1}, S_{4,1})$, $(R_1, S_{2,1}, S_{3,1})$ and $(R_1, S_{2,1}, S_{4,1})$.

From the action of $\a_B$, we see that the two first possibilities correspond to the two generators of the $\a_B$-invariant sublattice $A_1^*$ of $\MWl$. Hence, the orthogonal complement ${A_1^*}^{\oplus 2}$ is generated by the minimal sections with the four last possible intersection numbers. For each of these case, $\cp_2(\x)$ intersects $(R_0, S_{1,1}, S_{2,1}, S_{3,1}, S_{4,1})$; hence $\cp_2(\x) = \h$. That is, for any section corresponding to a minimal point of $\ker(\PH_2)$, we have that $\cp_2(\x) = \h$ and $d=2$.

Therefore, we obtain that $\ker(\PH_2)$ splits into the disjoint union of the 45 degree rotated sublattice $\ker(\PH_2)_{d=1} = {\langle 1 \rangle}^{\oplus 2} \oplus \IZ_2$ and its non-trivial coset $\ker(\PH_2)_{d=2}$.\\

\noindent\emph{Case 27.} $m=2$. $\ker(\PH_2) = MW = A_1^* \oplus \IZ_4$. $f_0$ is of type $I_2$, and $T = {A_3}^{\oplus 2} \oplus A_1$. Let $R_i$, $i=0,1$ be the components of $f_0$, and $S_{j,k}$, $j=1,2$ and $k=0,1,2,3$ be the components of the two other singular fibers with root lattices $A_3$.

The contributions for $f_0$ are either $0$ for the neutral component, or $1/2$ otherwise. For the two other fibers, the contributions are either $0$ for the neutral component, $3/4$ for the components $1$ and $3$, and $1$ for the middle component $2$.

There are three non-zero torsion sections of order $4$, which we denote by $\h$, $\h \boxplus \h$ and $\h \boxplus \h \boxplus \h$. From proposition \ref{p:torsion}, we get that $\h$ intersects $R_1$, $S_{1,1}$ and $S_{2,1}$, $\h \boxplus \h$ intersects $R_0$, $S_{1,2}$ and $S_{2,2}$, and $\h \boxplus \h \boxplus \h$ intersects $R_1$, $S_{1,3}$ and $S_{2,3}$.

$\a_B$ sends each component of $f_0$ to itself, and permutes the other components. That is,
$$
\a_B:~~\{ R_i, S_{1,k}, S_{2,k} \} \mapsto \{ R_i, S_{2,k}, S_{1,k} \},
$$
for $i=0,1$ and $k=0,1,2,3$. Hence, for any section $\x \in \ker(\PH_2)$, we know that $\cp_2(\x)$ must intersect the neutral component $R_0$. Thus, either $\cp_2(\x) = \s$, or $\cp_2(\x) = \h \boxplus \h$. In the latter case, $d=2$, since $\h \boxplus \h$ generates the subgroup $\IZ_2 \subset \IZ_4$. 

Now, take a section $\x \in \ker(\PH_2)$. It is easy to show that if $\x$ satisfies $\cp_2(\x) = \s$, then $\x \boxplus \h$ satisfies $\cp_2(\x \boxplus \h) = \h \boxplus \h$, $\x \boxplus \h \boxplus \h$ satisfies $\cp_2(\x \boxplus \h \boxplus \h) = \s$, and $\cp_2(\x \boxplus \h \boxplus \h \boxplus \h) = \h \boxplus \h$. Hence, we get that $AS = A_1^* \oplus \IZ_4$ splits into the disjoint union of the index 2 sublattice $\ker(\PH_2)_{d=1} = A_1^* \oplus \IZ_2$, and its non-trivial coset $\ker(\PH_2)_{d=2}$.
\end{proof}

\subsection{Free action on $f_\infty$}
We have now computed explicitly $\ker(\PH_m)$ for each cases in table \ref{t:pairs}, and found out when $d>1$ is possible. However, this is not the end of the story; we also want $\bra \t_B \ket = \bra t_\x \circ \a_B \ket$ to act freely on $f_\infty$.

\begin{proposition}\label{p:free}
Let $(B,\a_B)$ be a suitable $\s$-pair with a smooth $f_\infty$. Take a section $\x \in \ker(\PH_m)$, and define
the automorphism $\t_B = t_\x \circ \a_B$, of order $n$, with $n = d m$. Then the following statements are equivalent:
\begin{enumerate}
\item $\langle \t_B \rangle$ acts freely on $f_\infty$;
\item $\cp_i(\x)\big|_{\infty} \neq 0$, for $i=1, \ldots, n-1$;
\item $\xi \big|_{\infty}$ is a torsion point of the smooth elliptic curve of order precisely $n$.
\end{enumerate}
\end{proposition}

\begin{proof}
We know that $f_\infty$ is smooth, and that $\a_B$ fixes $f_\infty$ pointwise. Since translation by a non-zero section acts freely on a smooth elliptic curve, $\langle \t_B \rangle$ will be free on $f_\infty$ if and only if $\cp_i(\x)$ is non-zero on $f_\infty$ for $i=1,\ldots,n-1$. Now since $\a_B$ fixes $f_\infty$ pointwise, this can happen if and only if $\xi$ intersects $f_\infty$ at a torsion point of order precisely $n$.
\end{proof}


\subsection{Summary}

\label{s:summproc}

Hence, taking into account lemma \ref{l:perps} and proposition \ref{p:free}, we now know how to construct the set of allowed sections, which we call $AS$, for each case in table \ref{t:pairs}. Let us summarize how we proceed.

For each suitable $\s$-pair $(B,\a_B)$ listed in table \ref{t:pairs}, we first extract $\ker(\PH_m)$ using lemmas \ref{l:all} and \ref{l:perps}. Then, for the cases listed in table \ref{t:d} where $d > 1$ is possible, we split $\ker(\PH_m)$ into the disjoint union of $\ker(\PH_m)_{d=1}$ and $\ker(\PH_m)_{d > 1}$ according to the result of lemma \ref{l:d}. This gives us the set of sections which yield automorphisms $\t_B$ of order $n$ for each case.

Then we must remove from these sets the sections which do not intersect $f_\infty$ at a torsion point of order precisely $n$, according to proposition \ref{p:free}. This can be done as follows. Consider the map $(1-\alpha_B): MW \to MW$. Its image $\Im(1-\alpha_B)$ is a finite index subgroup of $\ker(\Phi_m)$. Since $\alpha_B$ fixes $f_\infty$ pointwise, any section in $\Im(1-\alpha_B)$ will necessarily intersect $f_\infty$ at zero.\footnote{Note also that when $d>1$ is allowed, $\Im(1-\alpha_B)$ is always a subgroup of $\ker(\PH_m)_{d=1}$, since for any 
$\zeta \in MW$, $\cp_m \left( (1-\alpha_B) \zeta \right) = \sigma$.}

As a consequence, all sections in a given coset of $\Im(1-\alpha_B)$ in $\ker(\Phi_m)$ will intersect $f_\infty$ at the same torsion point (or $0$).
This tells us that to obtain the set of allowed sections $AS$ from $\ker(\PH_m)$ for a given $n$ (hence a given $d$), we first need to compute the subgroup $\Im(1-\alpha_B) \subset \ker(\Phi_m)$, and then keep only the cosets containing sections intersecting $f_\infty$ at a torsion point of order precisely $n$. We thus obtain a complete list of pairs $(B,\t_B)$, with smooth $f_\infty$ and $\t_B$ acting freely on $f_\infty$.

\begin{remark}
Note that one can show that all sections in a given coset of $\Im(1-\alpha_B) \subset \ker(\Phi_m)$ lead to isomorphic quotients $B / \bra \tau_B \ket$, hence the number of non-isomorphic quotients, for a given $\alpha_B$, is at most the cardinality of the quotient group $\ker(\Phi_m) / \Im(1-\alpha_B)$. To see that, consider two automorphisms $\tau_B = t_{\xi} \circ \alpha_B$ and $\tau'_B = t_{\xi'} \circ \alpha_B$, with $\xi, \xi' \in \ker(\Phi_m)$. A sufficient condition for the quotients to be isomorphic is that there is an isomorphism $a:B \to B$ intertwining the $\tau$'s: $a \circ \tau_B = \tau_B' \circ a$. Take $a=t_{\zeta}$ for some section $\zeta$. Explicitly, the
condition then says:
\begin{equation}
t_{\zeta} \circ t_{\xi} \circ \alpha_B = t_{\xi'} \circ \alpha_B \circ t_{\zeta}.
\end{equation}
Since these two automorphisms have the same linearization, they
agree if the two of them take the zero section to the same section.
We thus get the condition
\begin{equation}
\xi' - \xi = (1-\alpha_B)\zeta,
\end{equation}
which implies that all sections in a given coset of $\Im(1-\alpha_B) \subset \ker(\Phi_m)$ produce isomorphic quotients.
\end{remark}


\section{The list(s)}\label{s:class}

Combining the results of the previous sections, we now have two lists:

\begin{enumerate}
\item A list of rational elliptic surfaces $B$ with automorphism groups generated by translation by torsion sections --- this was extracted directly from the torsion groups of the rational elliptic surfaces in Persson's list \cite{Miranda:1990, Persson:1990};
\item A list of pairs $(B,\t_B)$, where $\t_B$ is an automorphism of the second kind, $f_\infty$ is smooth and $\t_B$ acts freely on $f_\infty$. In other words, we have a list of rational elliptic surfaces $B$, with smooth $f_\infty$, and cyclic automorphism groups generated by automorphisms of the second kind acting freely on $f_\infty$.
\end{enumerate}

By combining appropriately these two lists, we are now in a position to produce a complete list of rational elliptic surfaces with finite automorphism groups (not necessarily cyclic) that can be lifted to free automorphism groups on smooth fiber products. The list naturally splits into two sublists, depending on whether the automorphism group acts trivially on the $\IP^1$ base or not.

\subsection{Trivial action on $\IP^1$}

In this case we obtain the following result.

\begin{proposition}
The list of rational elliptic surfaces with finite automorphism groups acting trivially on $\IP^1$ is given in table \ref{t:listtriv}. In this table we write down the automorphism group $G_B$ (which consists in translation by  torsion sections), the dimension ${\rm dim}$ of the moduli space, and the configuration of singular fibers at a generic point in the moduli space. We also write down the root lattice $T$ associated to the singular fibers and the Mordell-Weil group $MW$, for the generic configuration of singular fibers. 

It is understood that for cases $5$--$8$ there are various
specializations in the moduli space where some of the singular fibers
collide to produce different configurations of singular fibers, and that these specializations do not
necessarily have the same root lattice $T$ and Mordell-Weil group $MW$ as the
generic configuration. The specializations that have the same torsion group $\MWt$ as the generic configuration are shown in table \ref{t:listtrivs}. Additional degenerations are possible in which the
torsion group $\MWt$ becomes larger --- for example, both the $(\IZ_2)^2$ locus
and the $\IZ_4$ locus can degenerate to the $\IZ_4 \times \IZ_2$ locus. We did not list such degenerations under the subgroup, as they can readily be found listed
under the larger group.
\end{proposition}

\begin{table}[!htb]
\begin{center}
\begin{tabular}{|c|ccccc|}
\hline
$\#$&$G_B$ &${\rm dim}$&${\rm Singular~fibers}$&$T$&$MW$\\
\hline
$1$&$\IZ_3 \times \IZ_3$&$0$&$\{ 4 I_3 \}$&$A_2^{\oplus 4}$&$(\IZ_3)^2$\\
$2$&$\IZ_4 \times \IZ_2$&$0$&$\{ 2 I_4, 2 I_2 \}$&$(A_3 \oplus A_1)^{\oplus 2}$&$\IZ_4 \times \IZ_2$\\
$3$&$\IZ_6$&$0$&$\{I_6, I_3, I_2, I_1 \}$&$A_5 \oplus A_2 \oplus A_1$&$\IZ_6$\\
$4$&$\IZ_5$&$0$&$\{ 2 I_5, 2I_1 \}$&$A_4^{\oplus 2}$&$\IZ_5$\\
$5$&$\IZ_4$&$1$&$\{2 I_4, I_2, 2I_1 \}$&$A_3^{\oplus 2} \oplus A_1$&$A_1^* \oplus \IZ_4$\\
$6$&$\IZ_2 \times \IZ_2$&$2$&$\{ 6 I_2 \}$&$A_1^{\oplus 6}$&${A_1^*}^{\oplus 2} \oplus (\IZ_2)^2$\\
$7$&$\IZ_3$&$2$&$\{3 I_3, 3I_1 \}$&$A_2^{\oplus 3}$&$A_2^* \oplus \IZ_3$\\
$8$&$\IZ_2$&$4$&$\{ 4 I_2, 4 I_1 \}$&$A_1^{\oplus 4}$&$D_4^* \oplus \IZ_2$\\
\hline
\end{tabular}
\vspace{0.5em}
\caption{List of rational elliptic surfaces with finite automorphism groups acting trivially on $\IP^1$.} \label{t:listtriv}
\end{center}
\end{table}

\begin{table}[!htb]
\begin{center}
\begin{footnotesize}
\begin{tabular}{|c|c|l|}
\hline
$\#$&$G_B$&Specializations\\
\hline
$5$&$\IZ_4$&$\{ 2 I_4, I_2, 2I_1 \}$, $\{I_8,I_2,2I_1 \}$, $\{I_1^*, I_4, I_1\}$\\
$6$&$\IZ_2 \times \IZ_2$&$\{ 6 I_2 \}$, $\{I_4, 4I_2 \}$, $\{I_0^*, 3 I_2 \}$, $\{ I_2^*, 2 I_2 \}$, $\{ 2 I_0^* \}$\\
$7$&$\IZ_3$&$\{ 3 I_3, 3I_1\}$, $\{3 I_3, I_2, I_1\}, \{IV, 2 I_3, 2 I_1\}, \{I_6, I_3, 3I_1 \}$,\\
&&$\{ IV, 2 I_3, I_2 \}, \{ I_6, IV, 2 I_1\}, \{I_9, 3 I_1\}$, $\{ 3 IV \}, \{ IV^*, I_3, I_1 \}$, $\{ IV^*, IV\}$\\
$8$&$\IZ_4$&$\{ 4 I_2, 4 I_1 \}$, $\{ 5 I_2, 2 I_1 \}, \{III, 3 I_2, 3 I_1\}, \{I_4, 2 I_2, 4 I_1 \}$,\\
&& $\{ III, 4 I_2, I_1\}, \{2 III, 2 I_2, 2 I_1 \}, \{ I_3, 4 I_2, I_1 \}, \{ I_4, 3 I_2, 2I_1\}, \{I_4, III, I_2, 3 I_1\}, \{2 I_4, 4 I_1\}$,\\&&$ \{I_6, I_2, 4 I_1 \}$,
$\{ 2 III, 3 I_2\}, \{ I_3, III, 3I_2 \}, \{I_4, III, 2 I_2, I_1\}, \{ I_4, 2 III, 2 I_1 \}, \{I_4, I_3, 2I_2, I_1\}$,\\
&&$ \{I_0^*, 2 I_2, 2I_1\}, \{I_6, 2 I_2, 2I_1\}, \{I_6, III, 3 I_1 \}, \{I_8, 4 I_1\}$, $\{ 4 III \}, \{ I_4, 2 III, I_2\},$\\&&$ \{I_4, I_3, III, I_2 \}, \{I_0^*, III, I_2, I_1\}, \{I_0^*, I_4, 2 I_1 \}, \{I_6, III, I_2, I_1\}, \{I_1^*, 2 I_2, I_1\},$\\
&&$ \{I_2^*, I_2, 2I_1\}$,
$\{I_0^*, 2 III \}, \{I_1^*, III, I_2\}, \{I_2^*, III, I_1\}, \{III^*, I_2, I_1 \}, \{I_4^*, 2 I_1\}$, $\{III^*, III\}$\\
\hline
\end{tabular}
\end{footnotesize}
\vspace{0.5em}
\caption{Specializations for table \ref{t:listtriv}.} \label{t:listtrivs}
\end{center}
\end{table}

\begin{proof}
As mentioned earlier, this list can be obtained directly from Persson's list of rational elliptic surfaces \cite{Miranda:1990, Persson:1990} by keeping all the surfaces with non-trivial torsion groups $\MWt$. (Alternatively, one could work with the main table in \cite{Oguiso:1991}; but note that
the torsion group in their case 70 should be $\IZ_4$, rather than $(\IZ_2)^2$ as
they claim. This is responsible for the first specialization, $\{I_8,I_2,2I_1\}$,
of the $\IZ_4$ surface, as listed in the first line of our Table 7.) 

What remains is the calculation of the dimension of the moduli space,
\emph{i.e.} the number of complex deformations. For this we present first a
sketch, for which we are grateful to the referee, and then a detailed and
explicit calculation.

It is well-known that semi-stable rational elliptic surfaces with fixed
configuration and Mordell-Weil group of rank $r$ have moduli space of dimension $r$. In cases
5--7, the configuration uniquely determines the torsion group $\MWt$ (see \cite{Oguiso:1991} --- this follows
from the height pairing in \cite{Shioda:1990}); so only case 8 remains where both
$\IZ_2$-torsion and torsion free are possible. 

To deal with case 8, start with case 5 in table \ref{t:pairs}; this surface has configuration of singular fibers $\{III^*, I_2, I_1 \}$, a 2-torsion section and no moduli (it is an
extremal rational elliptic surface). We want to get rid of the $III^*$ fiber. According to table \ref{t:ram}, this requires a base change of degree $4$ with ramification of order $4$ at the $III^*$ fiber, to obtain the desired configuration $\{ 4 I_2, 4 I_1 \}$. However, in this case generally we do not know whether the pullbacked surface has $\IZ_2$-torsion or not. But here the situation is special; if the ramification index at the $III^*$ fiber is $(2,2)$, then we obtain two $I_0^*$ fibers in the pullback, which can then be eliminated by a quadratic twist (which transforms ``star"-fibers into ``unstarred"-fibers, see for instance \cite{Miranda:1990}). We thus obtain the desired configuration of singular fibers with $\IZ_2$-torsion in $MW$, since quadratic twists preserve 2-torsion sections. We can write down this base change explicitly as follows. We normalize the fiber $III^*$ to sit at $\infty$. We can choose three pre-images for the base change --- let us choose $0$ and $\infty$ as pre-images for $\infty$, and $1$ as one pre-image for $0$. Such base changes take the form
$$
t \mapsto (at^3+bt^2+ct+d) (t-1)/t^2.
$$ 
We deduce that the moduli space is at least four-dimensional. On the other
hand, by the previous argument, the dimension of the moduli space can at most
equal the rank of $MW$, which is 4 here. Hence the moduli space in this case is four-dimensional as
claimed. (This last argument actually works for all cases 5--8, by considering appropriate base changes.)

These dimensions can also be obtained by explicit calculations. We start
with an argument which occurs in the proof of Theorem 7.1 below. This says
that the number of deformations of a threefold $X$ which is a quotient of a Schoen's threefold is given by
$$
h = h_B + h_{B'} + e,
$$
where $h_B$ and $h_{B'}$ are the number of deformations of the $G = G_B = G_{B'}$ action on the two rational elliptic surfaces $B$ and $B'$, and $e = 3$ when $m=1$, while $e=1$ otherwise. Now, since
$$
H^2(X,\IQ) = {H^2(B,\IQ)^G \oplus H^2(B',\IQ)^G \over H^2(\IP^1, \IQ)},
$$ 
where the $G$ superscript denotes the $G$-invariant part, we obtain that the number of $G$ deformations of $B$ is given by
$$
h_B = \dim H^2 ( B, \IQ)^G - {e+1 \over 2}. 
$$

To compute the dimension of the $G$-invariant part of $H^2(B,\IQ)$, we use the Lefschetz fixed-point theorem. We obtain the general formula:
$$
\dim H^2 ( B, \IQ)^G = {12 + \sum_{i=1}^{n-1} f_i \over  n}   - 2,
$$
where $n$ is the order of $G$ and the $f_i$, $i=1, \ldots, n-1$ are the numbers of fixed points of the non-zero elements of $G$. Thus, by counting the number of fixed points of the action on $B$ of each element in $G$, we obtain $\dim H^2 ( B, \IQ)^G$, hence the number of deformations $h_B$. This gives the numbers $\dim$ in the table.

Finally, we need to check that all the cases listed in Table 7 are indeed specializations of the generic configuration. For most cases this is straightforward. For the surfaces involving star fibers, a simple way to do so is to count the number of fixed points as above for these subfamilies to obtain their number of deformations $h_B$. Recall that whenever the
$j$-invariant of an elliptic fibration is non-constant, the complex structure is determined locally by the positions in the base of the singular fibers. So if $h_B$ is bigger than the number of the singular fibers minus $3 = \dim PGL(2)$, this means that the surface is in a subspace of a bigger moduli space, hence it is a specialization of a more generic configuration. Using this simple argument we can prove that all the configurations of table 7 are specializations of the generic configuration, for each $G_B$. Let us do a few cases explicitly.

For case $\# 5$ with $G_B = \IZ_4$, we want to show that the configuration $\{I_1^*, I_4, I_1\}$ is a degeneration of the generic configuration $\{2 I_4, I_2, 2 I_1 \}$. The order $4$ elements of $\IZ_4$ have one fixed point on the $I_1^*$ fiber and one fixed point on the $I_1$ fiber, while the order $2$ element has 3 fixed points on $I_1^*$ and one fixed point on $I_1$. Hence we get that $h_B = 1$; but this is an extremal surface, which, according to \cite{Miranda:1986}, is unique. Hence it must be a specialization of the generic configuration; that is, we showed that the fibers $I_4, I_2$ and $I_1$ collide to give an $I_1^*$ fiber.

For case $\# 6$ with $G_B = \IZ_2 \times \IZ_2$, we need to show that $3 I_2$ collide to give a $I_0^*$. For the configuration $\{I_0^*, 3 I_2\}$, each element of $G_B$ has 2 fixed points on $I_0^*$ and 2 fixed points on one of the $I_2$ fibers. Hence $h_B = 2$, which means that it is indeed a degeneration of the generic configuration $\{ 6 I_2 \}$. For this case we also need to show that $4 I_2$'s collide to give an $I_2^*$. For the configuration $\{I_2^*, 2I_2 \}$, two elements have 2 fixed points on $I_2^*$ and two fixed points on one of the two $I_2$ fibers, and one element has 4 fixed points on $I_2^*$ and acts freely on the $I_2$'s. Hence $h_B=2$ again, and the surface is a specialization of the generic configuration.

For case $\# 7$ with $G_B = \IZ_3$, we need to check that $2I_3$ and $2 I_1$ collide to give a $IV^*$ fiber. For the configuration $\{IV^*, I_3, I_1 \}$, each non-zero elements of $G_B$ has 2 fixed points on $IV^*$ and one fixed point on $I_1$. Hence $h_B = 2$, that is it is a specialization of the generic configuration with $\{3 I_3, 3 I_1 \}$.

Finally, for case $\# 8 $ with $G_B = \IZ_2$ it remains to be checked that $3I_2$ and $3 I_1$ collide to give a $III^*$ fiber, and that $4I_2$ and $2I_1$ collide to give a $I_4^*$ fiber. For the configuration $\{III^*, I_2, I_1\}$, the non-zero element has 3 fixed points on $III^*$ and one fixed point on $I_1$, hence $h_B = 4$ and it is a specialization of the generic configuration. For the configuarion $\{I_4^*, 2I_1\}$, the non-zero element has 2 fixed points on $I_4^*$ and one fixed point on each $I_1$'s, hence $h_B = 4$ and it is also a specialization of the generic configuration.
\end{proof}

\subsection{Non-trivial action on $\IP^1$}

In this case a little bit more work is needed. We obtained in the previous section a list of suitable pairs $(B,\t_B)$. However, some of these surfaces also have non-trivial torsion group $\MWt$. In this case, we can form ``mixed" finite automorphism groups, with one factor generated by an automorphism of the second kind and the other factors generated by translation by torsion sections. Taking this into account, we obtain the following lemma.

\begin{proposition}\label{p:list}
The list of rational elliptic surfaces with $f_\infty$ smooth and finite automorphism groups $G_B$ acting non-trivially on $\IP^1$ and freely on $f_\infty$ is given in tables \ref{t:list1} and \ref{t:list2}. For each entry in the tables, we write down, following the notation of the previous sections:
\begin{itemize}
\item the order $m$ of the linearization $\alpha_B$;
\item the order $d$ of the torsion section $\cp_m (\x)$;
\item the dimension ${\rm dim}$ of the moduli space of complex deformations preserving the automorphism group;
\item the generic configuration of singular fibers in the moduli space (the fiber in bold typeface corresponds to the special fiber $f_0$; hence if there is no bold singular fiber it means that $f_0$ is smooth);
\item the root lattice $T$ associated to the singular fibers of $B$;
\item the Mordell-Weil group $MW$ of $B$;
\item the $\a_B$-invariant subgroup $MW^{\a_B}$ of $MW$;
\item the subgroup of sections $\ker(\PH_m) \subset MW$;
\item when $d > 1$ is possible, the subgroup $\ker(\PH_m)_{d=1} \subset \ker(\PH_m)$ --- the set $\ker(\PH_m)_{d>1}$ is then the union of the non-trivial cosets of $\ker(\PH_m)_{d=1} \subset \ker(\PH_m)$;
\item the corresponding case number ${\rm num}$ of table \ref{t:pairs}.
\end{itemize}
 For each entry in the table it is understood that there may be various specializations in the moduli space where some of the singular fibers collide to produce different configurations of singular fibers, with the same associated root lattice $T$ and Mordell-Weil group $MW$. These specializations are shown in table \ref{t:lists}; we use semicolons to segregate strata of a given dimension, which are then separated by commas. In table \ref{t:lists} we only list the cases where there are specializations.
\end{proposition}

\begin{remark}
\label{r:as}
Recall that to get the set of allowed sections $AS$ in each case, we must keep only the cosets of $\Im(1-\alpha_B) \subset \ker(\PH_m)$ which contain sections intersecting $f_\infty$ at a torsion point of order precisely $n$, with $n=d m$. One can show that $AS$ is non-empty for all lines in tables \ref{t:list1} and \ref{t:list2}; we sketch a case-by-case proof in the Appendix.
\end{remark}

\begin{remark}
Note that, as explained in section \ref{s:simple}, there are codimension one subspaces in the moduli space where some of the singular fibers collide to give a singular fiber at $\infty$. However the full automorphism group may not act freely on $f_\infty$ at these points, depending on the singular type of $f_\infty$. In any case, we excluded these cases from our analysis for the reasons pointed out in section \ref{s:simple}, and considered only generic configurations of singular fibers in a given moduli space, where $f_\infty$ is smooth.
\end{remark}

\begin{remark}
Note that various entries in the tables correspond to specializations of other entries with the same automorphism group. We decided to write separate entries in tables \ref{t:list1} and \ref{t:list2} whenever the surfaces have different Mordell-Weil groups, since the set of allowed sections $AS$ in the construction of the automorphism $\t_B$ is then different, resulting in different
possibilities for spectral covers and eventually vector bundles, cf. \cite{Bouchard:tocome}. We will come back to this in remark \ref{r:specializations}.
\end{remark}

\begin{proof}
We obtained in the previous section a list of suitable $(B,\t_B)$, where $\t_B$ is an automorphism of the second kind, that is with $m>1$, $f_\infty$ is smooth, and $\t_B$ acts freely on $f_\infty$. For the cases where $\MWt$ is non-trivial, one can combine this cyclic automorphism group with the automorphism group generated by translation by torsion sections to obtain a (perhaps non-cyclic) automorphism group.

For the group to be abelian, the generators of the cyclic factors must commute, that is, the torsion sections must be $\a_B$-invariant (since translation by sections always commute). However, by looking at table \ref{t:pairs} it is easy to see that for any pair $(B,\a_B)$ in our list, the torsion sections are always $\a_B$-invariant; hence the automorphism groups that we obtain are always abelian.

We also want the full group to act freely on $f_\infty$. We know that the cyclic automorphism group generated by an automorphism of the second kind acts freely on $f_\infty$, by construction. Moreover, any translation by a non-zero torsion section also acts freely on the smooth fiber $f_\infty$. But we must make sure that no element of $G_B$ acts trivially on $f_\infty$. The only subtle cases in that respect are cases 22 and 26. Here it may seem that one could obtain a $(\IZ_2)^3$ automorphism group by combining the $\IZ_2$ automorpism group acting on $\IP^1$ with a torsion $(\IZ_2)^2$. However, there is no free $(\IZ_2)^3$ action on a smooth elliptic curve, hence at least one element of this group must fix the smooth fiber $f_\infty$ pointwise. That is, only a subgroup $(\IZ_2)^2 \subset (\IZ_2)^3$ acts freely on $f_\infty$.

To compute the dimension of the moduli space of complex deformations, we first follow the same argument as in the proof of proposition 6.1, but with $e=1$ since now $m > 1$. This gives the dimension of the moduli space for the generic configuration of singular fibers. For the entries in the table that correspond to specializations of other entries, we proceed as follows.

Assume we are given a family of elliptic surfaces, determined by a specific configuration of singular fibers, and that the positions of those singular fibers are subject to a number $c$ of constraints. Consider another such family which is in the boundary of the former, \emph{i.e.} it is given by a singular configuration which is a degeneration of the first. Then the number of constraints on the new singular fibers is $\leq c$. In particular, if $c=0$ for the big family, then $c=0$ for any stratum in it. This is the case for all the families in table 8, in the sense that the position of the singular fibers of the quotient surfaces are unrestricted. Hence we get the dimension of the moduli space of these subfamilies by counting the number of singular fibers of the quotient surface minus $3 = \dim PGL(2)$. This also accounts for the stratification of the configurations in table 9.

Taking these remarks into account the list follows from the results of the previous sections. Let us simply end the proof by commenting on a few cases:\\

\noindent\emph{Case 11}. In this case one may think that by choosing a section $\x \in \ker(\PH_2)_{d=2}$
to build $\t_B$, and combining with translation by a torsion section of order $4$, one would get an automorphism group $\IZ_4 \times \IZ_4$. However, this is not true since $\t_B^2 = t_{\h \boxplus \h}$, where $\h \boxplus \h$ is the order $2$ torsion section; thus only a $\IZ_4 \times \IZ_2$ group acts faithfully. This kind of analysis applies each time $d>1$ is allowed.
\\

\noindent\emph{Cases 11, 14 and 15}. In these cases we can realize the automorphism group either by taking a section $\x \in \ker(\PH_m)_{d=1}$ or by taking a section $\x \in \ker(\PH_m)_{d > 1}$, and combining with translation by a torsion section. But these two realizations of the automorphism group are in fact equivalent, since they only differ by a linear combination of the generators.\\

\noindent\emph{Cases 12, 19, 20 and respectively 26, 24, 28}. These are the same rational elliptic surfaces. The three first cases correspond to the automorphism groups obtained by choosing $\x \in \ker(\PH_2)_{d=2}$, while the three last cases correspond to the choice $\x \in \ker(\PH_2)_{d=1}$.\\

\noindent\emph{Cases 22, 26, 29}. Notice that in these cases one can also realize the $(\IZ_2)^2$ automorphism group only by translation by torsion sections ($m=1$).
\end{proof}

\newpage
${}$
\vfill

\begin{table}[!b]
\begin{center}
\begin{footnotesize}
\begin{tabular}{|c|cccccccccc|}
\hline\hline
\multicolumn{11}{|c|}{$\IZ_3 \times \IZ_3$}\\[0.2cm]
\hline
$\#$&$m$&$d$&${\rm dim}$&${\rm Sing.~fibers}$&$T$&$MW$&$MW^{\a_B}$&$\ker(\PH_m)$&$\ker(\PH_m)_{d=1}$&$num$\\
\hline
9&$3$&$1$&$1$&$\{3 I_3, 3I_1\}$&${A_2}^{\oplus 3}$&$A_2^* \oplus \IZ_3$&$\IZ_3$&$MW$&-&$11$\\
\hline\hline
\multicolumn{11}{|c|}{$\IZ_4 \times \IZ_2$}\\[0.2cm]
\hline
10&$4$&$1$&$1$&$\{4I_2, 4I_1 \}$&${A_1}^{\oplus 4}$&$D_4^* \oplus \IZ_2$&$\IZ_2$&$MW$&-&$5$\\

11&$2$&$2/1$&$1$&$\{ \mathbf{I_2}, 2 I_4, 2 I_1 \}$&${A_3}^{\oplus 2} \oplus A_1$&$A_1^* \oplus \IZ_4$&$\IZ_4$&$MW$&$A_1^* \oplus \IZ_2$&$27$\\

12&$2$&$2$&$1$&$\{ {\bf I_4}, 4 I_2 \}$&$A_3 \oplus {A_1}^{\oplus 4}$&$\langle 1/4 \rangle \oplus (\IZ_2)^2$&$(\IZ_2)^2$&$MW$&$\langle 1 \rangle \oplus (\IZ_2)^2$&$22$\\

\hline\hline
\multicolumn{11}{|c|}{$\IZ_6$}\\[0.2cm]
\hline
13&$6$&$1$&$1$&$\{ 12 I_1 \}$&$0$&$E_8$&$0$&$MW$&-&$1$\\

14&$3$&$2/1$&$1$&$\{ {\bf III}, 3 I_2, 3I_1 \}$&${A_1}^{\oplus 4}$&$D_4^* \oplus \IZ_2$&$\IZ_2$&$MW$&$D_4^*$&$8$\\

15&$2$&$3/1$&$1$&$\{ {\bf IV}, 2 I_3, 2I_1 \}$&${A_2}^{\oplus 3}$&$A_2^* \oplus \IZ_3$&$\IZ_3$&$MW$&$A_2^*$&$17$\\

\hline\hline
\multicolumn{11}{|c|}{$\IZ_5$}\\[0.2cm]
\hline
16&$5$&$1$&$1$&$\{ {\bf II}, 10 I_1 \}$&$0$&$E_8$&$0$&$MW$&-&$2$\\

\hline\hline
\multicolumn{11}{|c|}{$\IZ_4$}\\[0.2cm]
\hline

17&$4$&$1$&$2$&$\{ 12 I_1 \}$&$0$&$E_8$&$A_1$&$E_7$&-&$4$\\

18&$4$&$1$&$1$&$\{ {\bf IV}, 8 I_1 \}$&$A_2$&$E_6^*$&$0$&$MW$&-&$3$\\

19&$2$&$2$&$2$&$\{ {\bf I_2}, 4 I_2, 2 I_1 \}$&${A_1}^{\oplus 5}$&${A_1^*}^{\oplus 3} \oplus \IZ_2$&$A_1^* \oplus \IZ_2$&${A_1^*}^{\oplus 2} \oplus \IZ_2$&$\langle 1 \rangle^{\oplus 2} \oplus \IZ_2$&$26$\\

20&$2$&$2$&$1$&$\{ {\bf I_0^*}, 2I_2, 2I_1 \}$&$D_4 \oplus {A_1}^{\oplus 2}$&${A_1^*}^{\oplus 2} \oplus \IZ_2$&$\IZ_2$&$MW$&$\langle 1 \rangle^{\oplus 2} \oplus \IZ_2$&$14$\\

\hline\hline
\multicolumn{11}{|c|}{$\IZ_2 \times \IZ_2$}\\[0.2cm]
\hline

21&$2$&$1$&$3$&$\{ 4 I_2, 4I_1 \}$&${A_1}^{\oplus 4}$&$D_4^* \oplus \IZ_2$&$\langle 1 \rangle^{\oplus 2} \oplus \IZ_2$&$\langle 1 \rangle^{\oplus 2} \oplus \IZ_2$&-&$31$\\

22&$2$&$1$&$2$&$\{ 6 I_2 \}$&${A_1}^{\oplus 6}$&${A_1^*}^{\oplus 2} \oplus (\IZ_2)^2$&$\langle 1 \rangle \oplus (\IZ_2)^2$&$\langle 1 \rangle \oplus (\IZ_2)^2$&-&$32$\\

23&$2$&$1$&$2$&$\{ 2 I_4, 4 I_1 \}$&${A_3}^{\oplus 2}$&${A_1^*}^{\oplus 2} \oplus \IZ_2$&$A_1^* \oplus \IZ_2$&$A_1^* \oplus \IZ_2$&-&$33$\\

24&$2$&$1$&$2$&$\{ {\bf I_2}, 4 I_2, 2 I_1 \}$&${A_1}^{\oplus 5}$&${A_1^*}^{\oplus 3} \oplus \IZ_2$&$A_1^* \oplus \IZ_2$&${A_1^*}^{\oplus 2} \oplus \IZ_2$&$\langle 1 \rangle^{\oplus 2} \oplus \IZ_2$&$26$\\

25&$2$&$1$&$2$&$\{ {\bf I_4}, 2 I_2, 4I_1 \}$&$A_3 \oplus {A_1}^{\oplus 2}$&$A_3^* \oplus \IZ_2$&$\langle 1 \rangle \oplus \IZ_2$&$\langle 1 \rangle^{\oplus 2} \oplus \IZ_2$&-&$21$\\

26&$2$&$1$&$1$&$\{ {\bf I_4}, 4 I_2 \}$&$A_3 \oplus {A_1}^{\oplus 4}$&$\langle 1/4 \rangle \oplus (\IZ_2)^2$&$(\IZ_2)^2$&$MW$&$\langle 1 \rangle \oplus (\IZ_2)^2$&$22$\\

27&$2$&$1$&$1$&$\{ {\bf I_8}, 4 I_1 \}$&$A_7$&$A_1^* \oplus \IZ_2$&$\IZ_2$&$MW$&-&$18$\\

28&$2$&$1$&$1$&$\{ {\bf I_0^*}, 2I_2, 2I_1 \}$&$D_4 \oplus {A_1}^{\oplus 2}$&${A_1^*}^{\oplus 2} \oplus \IZ_2$&$\IZ_2$&$MW$&$\langle 1 \rangle^{\oplus 2} \oplus \IZ_2$&$14$\\

29&$2$&$1$&$0$&$\{ 2 I_0^* \}$&${D_4}^{\oplus 2}$&$(\IZ_2)^2$&$(\IZ_2)^2$&$MW$&-&$34$\\

\hline\hline
\multicolumn{11}{|c|}{$\IZ_3$}\\[0.2cm]
\hline

30&$3$&$1$&$3$&$\{ 12 I_1 \}$&$0$&$E_8$&$A_2$&$E_6$&-&$9$\\

31&$3$&$1$&$2$&$\{ 3 I_2, 6 I_1 \}$&${A_1}^{\oplus 3}$&$D_4^* \oplus A_1^*$&$A_1^*$&$D_4^*$&-&$10$\\

32&$3$&$1$&$2$&$\{ {\bf III}, 9 I_1 \}$&$A_1$&$E_7^*$&$\langle 3/2 \rangle$&$E_6^*$&-&$7$\\

33&$3$&$1$&$1$&$\{ {\bf I_0^*}, 6 I_1 \}$&$D_4$&$D_4^*$&$0$&$MW$&-&$6$\\

\hline\hline
\end{tabular}
\end{footnotesize}
\vspace{0.5em}
\caption{List of rational elliptic surfaces with finite automorphism groups satisfying the criteria of proposition \ref{p:list}.}\label{t:list1}
\end{center}
\end{table}

\newpage

\begin{table}[!t]
\begin{center}
\begin{footnotesize}
\begin{tabular}{|c|cccccccccc|}
\hline\hline
\multicolumn{11}{|c|}{$\IZ_2$}\\[0.2cm]
\hline
$\#$&$m$&$d$&${\rm dim}$&${\rm Sing.~fibers}$&$T$&$MW$&$MW^{\a_B}$&$\ker(\PH_m)$&$\ker(\PH_m)_{d=1}$&$num$\\
\hline

34&$2$&$1$&$5$&$\{ 12 I_1 \}$&$0$&$E_8$&$D_4$&$D_4$&-&$28$\\

35&$2$&$1$&$4$&$\{2 I_2, 8 I_1 \}$&${A_1}^{\oplus 2}$&$D_6^*$&${\langle 1 \rangle}^{\oplus 3}$&$\langle 1 \rangle^{\oplus 3}$&-&$29$\\

36&$2$&$1$&$4$&$\{ {\bf I_2}, 10 I_1 \}$&$A_1$&$E_7^*$&$2 A_3^*$&$D_4$&-&$23$\\

37&$2$&$1$&$3$&$\{ 2 I_3, 6 I_1 \}$&${A_2}^{\oplus 2}$&${A_2^*}^{\oplus 2}$&$A_2^*$&$A_2^*$&-&$30$\\

38&$2$&$1$&$3$&$\{ {\bf I_2}, 2 I_2, 6 I_1 \}$&${A_1}^{\oplus 3}$&$D_4^* \oplus A_1^*$&$A_1^* \oplus \langle 1 \rangle$&$\langle 1 \rangle^{\oplus 3}$&-&$24$\\

39&$2$&$1$&$3$&$\{ {\bf I_4}, 8 I_1 \}$&$A_3$&$D_5^*$&$\langle 1 \rangle^{\oplus 2}$&$\langle 1 \rangle^{\oplus 3}$&-&$20$\\

40&$2$&$1$&$3$&$\{ {\bf IV}, 8 I_1 \}$&$A_2$&$E_6^*$&$2 A_2^*$&$D_4$&-&$15$\\

41&$2$&$1$&$2$&$\{ {\bf I_2}, 2 I_3, 4I_1 \}$&${A_2}^{\oplus 2} \oplus A_1 $&$A_2^* \oplus \langle 1/6 \rangle$&$\langle 1/6 \rangle$&$A_2^*$&-&$25$\\

42&$2$&$1$&$2$&$\{ {\bf I_6}, 6 I_1 \}$&$A_5$&$A_2^* \oplus A_1^*$&$A_1^*$&$A_2^*$&-&$19$\\

43&$2$&$1$&$2$&$\{ {\bf IV}, 2 I_2, 4I_1 \}$&$A_2 \oplus {A_1}^{\oplus 2}$&$\frac{1}{6}\begin{pmatrix}2 & 1 &0&-1\\1&5&3&1\\0&3&6&3\\-1&1&3&5 \end{pmatrix}$&$\langle 1/3 \rangle$&$\langle 1 \rangle^{\oplus 3}$&-&$16$\\

44&$2$&$1$&$2$&$\{ {\bf I_0^*}, 6 I_1 \}$&$D_4$&$D_4^*$&$\langle 1 \rangle$&$\langle 1 \rangle^{\oplus 3}$&-&$13$\\

45&$2$&$1$&$1$&$\{ {\bf IV^*}, 4 I_1 \}$&$E_6$&$A_2^*$&$0$&$MW$&-&$12$\\

\hline\hline
\end{tabular}
\end{footnotesize}
\vspace{0.5em}
\caption{List of rational elliptic surfaces with finite automorphism groups satisfying the criteria of proposition \ref{p:list} (continued).}\label{t:list2}
\end{center}
\end{table}

\begin{table}[!hbt]
\begin{center}
\begin{footnotesize}
\begin{tabular}{|c|c|l|}
\hline
$\#$&$G_B$&Specializations\\
\hline
$9$&$\IZ_3 \times \IZ_3$&$\{3 I_3, 3 I_1\}$; $\{3 IV \}$\\
$10$&$\IZ_4 \times \IZ_2$&$\{4 I_2, 4 I_1 \}$; $\{ 4 III \}$\\
$13$&$\IZ_6$&$\{ 12 I_1 \}$; $\{ 6 II \}$\\
$14$&&$\{III, 3 I_2, 3 I_1 \}$; $\{4 III \}$\\
$15$&&$\{IV, 2 I_3, 2 I_1 \}$; $\{3 IV \}$\\
$16$&$\IZ_5$&$\{II, 10 I_1 \}$; $\{6 II\}$\\
$17$&$\IZ_4$&$\{12 I_1 \}$; $\{4II, 4I_1\}$\\
$18$&&$\{IV, 8 I_1\}$; $\{IV, 4 II\}$\\
$19$&&$\{ 5 I_2, 2 I_1\}$; $\{ 3 I_2, 2 III\}$\\
$20$&&$\{ I_0^*, 2 I_2, 2 I_1\}$; $\{ I_0^*, 2 III\}$\\
$21$&$\IZ_2 \times \IZ_2$&$\{4 I_2, 4 I_1\}$; $\{ 2 III, 2 I_2, 2 I_1\}$; $\{ 4 III\}$\\
$24$&&$\{ 5 I_2, 2 I_1\}$; $\{ 3 I_2, 2 III\}$\\
$25$&&$\{I_4, 2 I_2, 4 I_1\}$; $\{I_4, 2 III, 2 I_1\}$\\
$28$&&$\{ I_0^*, 2 I_2, 2 I_1\}$; $\{ I_0^*, 2 III\}$\\
$30$&$\IZ_3$&$\{ 12 I_1 \}$; $\{3 II, 6 I_1\}$; $\{ 6 II\}$\\
$31$&&$\{3 I_2, 6 I_1 \}$; $\{3 I_2, 3 II\}, \{ 3III, 3 I_1\}$\\
$32$&&$\{III, 9 I_1\}$; $\{III, 3 II, 3 I_1\}$\\
$33$&&$\{I_0^*, 6 I_1 \}$; $\{I_0^*, 3 II\}$\\
$34$&$\IZ_2$&$\{12 I_1\}$; $\{ 2 II, 8 I_1\}$; $\{4II, 4 I_1\}$; $\{6 II\}$\\
$35$&&$\{2 I_2, 8 I_1\}$; $\{2 I_2, 2 II, 4 I_1\}, \{2 III, 6 I_1\}$; $\{2 I_2, 4 II\},\{2 III, 2 II, 2 I_1\}$\\
$36$&&$\{I_2, 10 I_1\}$; $\{I_2, 2 II, 6 I_1\}$; $\{I_2, 4 II, 2 I_1\}$\\
$37$&&$\{2 I_3, 6 I_1\}$; $\{2 I_3, 2 II, 2 I_1\}, \{2 IV, 4 I_1\}$; $\{2 IV, 2 II\}$\\
$38$&&$\{3 I_2, 6 I_1\}$; $\{3 I_2, 2 II, 2 I_1\}, \{I_2, 2 III, 4 I_1\}$; $\{I_2, 2 III, 2 II\}$\\
$39$&&$\{I_4, 8 I_1\}$; $\{I_4, 2 II, 4 I_1\}$; $\{ I_4, 4 II\}$\\
$40$&&$\{IV, 8 I_1\}$; $\{IV, 2 II, 4 I_1\}$; $\{IV, 4 II\}$\\
$41$&&$\{2 I_3, I_2, 4 I_1\}$; $\{2 I_3, I_2, 2 II\}, \{2 IV, I_2, 2 I_1\}$\\
$42$&&$\{I_6, 6 I_1\}$; $\{I_6, 2 II, 2 I_1\}$\\
$43$&&$\{IV, 2 I_2, 4 I_1\}$; $\{IV, 2 I_2, 2 II\}, \{ IV, 2 III, 2 I_1\}$\\
$44$&&$\{I_0^*, 6 I_1\}$; $\{I_0^*, 2 II, 2 I_1\}$\\
$45$&&$\{IV^*, 4 I_1\}$; $\{IV^*, 2 II\}$\\
\hline
\end{tabular}
\end{footnotesize}
\vspace{0.5em}
\caption{Specializations for tables \ref{t:list1} and \ref{t:list2}.} \label{t:lists}
\end{center}
\end{table}

\begin{remark}\label{r:specializations}
By looking at the configurations of singular fibers in tables \ref{t:list1} and \ref{t:list2}, it is straightforward to see that:
\begin{itemize}
\item case 18 is a specialization of case 17;
\item case 20 is a specialization of case 19;
\item cases 22--29 are specializations of case 21;
\item cases 31--33 are specializations of case 30;
\item cases 35--45 are specializations of case 34.
\end{itemize}
\end{remark}

Tables \ref{t:listtriv}--\ref{t:lists} provide a complete list of rational elliptic surfaces with finite automorphism groups that can be lifted to free automorphism groups on smooth fiber products. Then, given any rational elliptic surface $B$ in the tables above, we can quotient by the finite automorphism group (or a subgroup thereof) to obtain a new rational elliptic surface. 

\begin{remark}
For the cases with trivial action on $\IP^1$ given in table \ref{t:listtriv}, it is interesting to note that, for a surface with the generic configuration of singular fibers, the quotient surface --- by the whole torsion group --- also has the same configuration of singular fibers.
\end{remark}

\section{Non-simply connected Calabi-Yau threefolds}\label{s:examples}

\subsection{The main theorem}
Finally, we apply our results on automorphisms of rational elliptic
surfaces to the classification of the corresponding non-simply connected Calabi-Yau
threefolds.

\begin{theorem}\label{t:main}
The moduli space of quotients of smooth fiber products
$\tX := B  \times_{\IP^1} B'$ (Schoen threefolds) by free, fiber-preserving finite
group actions has components as indicated in table \ref{t:threefolds}.

In the table, $G$ is the group acting freely on $\tX$, $m$ is the order of the
image of $G$ in ${\rm Aut}(\IP^1)$, and $h$ is the Hodge number $h:=h^{1,1}=h^{2,1}$ of the
quotient threefold $X:= \tX / G$. The cases refer to tables \ref{t:listtriv}, \ref{t:list1} and \ref{t:list2}. When
we list a case as $n_1 \times n_2$ we refer to the fiber product of $B$ of case $n_1$ (in
the tables \ref{t:listtriv}--\ref{t:list2}) with $B'$ of case $n_2$.  When we list a case as $n$, it
means $n \times n$. When several cases are possible, we use semicolons to segregate
strata of a given dimension, which are then separated by commas.
\end{theorem}

\begin{table}[!htb]
\begin{center}
\begin{tabular}{|c|ccc|}
\hline
$G$&$m$ &$h$&${\rm Cases}$\\
\hline\hline
$\IZ_3 \times \IZ_3$&$3$&$3$&$9$\\
&$1$&$3$&$1$\\
\hline
$\IZ_4 \times \IZ_2$&$4$&$3$&$10$\\
&$2$&$3$&$11$\\
&$2$&$3$&$12$\\
&$2$&$3$&$11 \times 12$\\
&$1$&$3$&$2$\\
\hline
$\IZ_6$&$6$&$3$&$13$\\
&$3$&$3$&$14$\\
&$2$&$3$&$15$\\
&$1$&$3$&$3$\\
\hline
$\IZ_5$&$5$&$3$&$16$\\
&$1$&$3$&$4$\\
\hline
$\IZ_4$&$4$&$5$&$17;10,18$\\
&$2$&$5$&$19;11,12,20$\\
&$1$&$5$&$5$\\
\hline
$\IZ_2 \times \IZ_2$&$2$&$7$&$21;22-25;11,15,26-28;29$\\
&$1$&$7$&$6$\\
\hline
$\IZ_3$&$3$&$7$&$30;31,32;9,14,33$\\
&$1$&$7$&$7$\\
\hline
$\IZ_2$&$2$&$11$&$34;35,36;21,37-40;22-25,41-44;11,15,26-28,45;29$\\
&$1$&$11$&$8$\\
\hline
\end{tabular}
\vspace{0.5em}
\caption{Classification of non-simply connected Calabi-Yau threefolds constructed as free, fiber-preserving quotients of smooth fiber products $B  \times_{\IP^1} B'$.} \label{t:threefolds}
\end{center}
\end{table}

\begin{proof}
By assumption, $X= \tX/G$ is the quotient by some finite group $G$ of a smooth
Schoen threefold $\tX =B \times_{\IP^1} B'$, where $B$ and $B'$ must be listed in tables \ref{t:listtriv}--\ref{t:list2}, with the same $G$ and the same $m$. What remains is to
determine when two such $X$'s live in the same moduli space, and to
calculate the Hodge numbers. 

For the first task, we simply need to
determine which configurations in tables \ref{t:list1} and \ref{t:list2} can specialize to which other
configurations. This was presented in remark \ref{r:specializations}. 

As to the Hodge numbers, the smooth Schoen
threefolds $\tX$ have $h^{1,1}(\tX)=h^{2,1}(\tX)=19$, so their Euler
characteristic vanishes. This latter property descends to $X$, so
$h^{1,1}(X)=h^{2,1}(X):=h$. We determine this by counting complex moduli.
Any deformation of $X$ lifts to a deformation of its universal cover $\tX$,
which in turn is given by a deformation of the surfaces $B$ and $B'$. Now the
equivalence relation on $\tX$ which gives the quotient $X$ is etale, in fact it
is isomorphic to the product $\tX \times G$. This property is clearly invariant
under deformations, so the deformed $X$ must be the quotient of a deformed
$\tX$ by a free action of the same group $G$. It must therefore be included in
tables \ref{t:listtriv}--\ref{t:list2}. It follows that $h$, the number of complex moduli of
our Calabi-Yau $X$, is given as:
$$
h = h_B + h_{B'} + e
$$
where $h_B$, $h_{B'}$ are the numbers of deformations of the $G$ action on $B$, $B'$
respectively, which are indicated in tables \ref{t:listtriv}--\ref{t:list2}, and $e$ gives the
number of extra parameters coming from an identification of the $\IP^1$ bases
of $B$ and $B'$:
\begin{itemize}
\item $e=\dim(PGL(2))=3$ in case $m=1$;
\item $e=1$ if $m>1$ (in this case we lose $2$ of the $3$ parameters since the points $0,\infty$ must go to each other).
\end{itemize}
\end{proof}

\subsection{Examples}

In this subsection we give a few examples of non-simply connected Calabi-Yau threefolds that we constructed in this paper.

\subsubsection{Calabi-Yau threefolds with fundamental group $\p_1(X) \cong \IZ_3 \times \IZ_3$}

According to theorem \ref{t:main}, we constructed two Calabi-Yau threefolds with $\p_1(X)\cong \IZ_3 \times \IZ_3$.

The first one has Hodge numbers $h=3$, and correspond to the fiber product of two rational elliptic surfaces in the one-parameter family given by case $\# 9$ of table \ref{t:list1}. These rational elliptic surfaces have $T = A_2^{\oplus 3}$, with smooth $f_0$; one $\IZ_3$ is generated by translation by a torsion section, and the other $\IZ_3$ acts faithfully on the $\IP^1$ base. The generic configuration of singular fibers is $\{ 3I_3, 3I_1\}$. This is precisely the Calabi-Yau threefold studied in great detail in \cite{Braun:2004xv}.

The second Calabi-Yau threefold also has Hodge numbers $h=3$ and consists in the fiber product of two rational elliptic surfaces corresponding to case $\# 1$ of table \ref{t:listtriv}. Here, the full automorphism group acts trivially on $\IP^1$. It would be interesting to investigate the construction of standard-model bundles on this threefold.

\subsubsection{Calabi-Yau threefolds with fundamental group $\p_1(X) \cong \IZ_6$}

According to theorem \ref{t:main}, we get four Calabi-Yau threefolds with $\p_1(X) \cong \IZ_6$; all have Hodge numbers $h=3$. We are presently studying the construction of standard-model bundles on these manifolds and will report on it in further publication \cite{Bouchard:tocome}.

\subsubsection{Calabi-Yau threefolds with fundamental group $\p_1(X) \cong \IZ_2$}

According to theorem \ref{t:main} we get two Calabi-Yau threefolds with $\p_1(X) \cong \IZ_2$, with Hodge numbers $h=11$. Here we simply note that the first one is the Calabi-Yau threefold that was used in \cite{Bouchard:2005ag, Bouchard:2006dn, Donagi:2000fw, Donagi:2000zs} to construct standard-model bundles. More precisely, in \cite{Bouchard:2005ag, Bouchard:2006dn, Donagi:2000fw, Donagi:2000zs} the covering Calabi-Yau threefold $\tX$ was constructed as a smooth fiber product of two rational elliptic surfaces with configuration of singular fibers $2 I_2, 8 I_1$, with smooth $f_0$ and $\a_B \simeq \IZ_2$; this is the generic configuration in the four-parameter family given by case $\# 35$ in table \ref{t:list2}.

\subsubsection{Schoen's constructions}

In section 9 of \cite{Schoen:1988}, Schoen constructs four different
smooth fiber products with free $\IZ_n$ action, with $n=2,3,4,6$
respectively. He considers four types of rational elliptic surfaces,
and
then fiber products of two rational elliptic surfaces of the same
type.
The four types of rational elliptic surfaces he considers have smooth
fibers at $0$ and $\infty$, and admit $\IZ_n$ actions with fixed
points
only on the smooth fiber $f_0$. The $\IZ_n$ actions act faithfully on
the $\IP^1$ base, hence in our notation $m=n$, that is in each case
the
order of $\a_B$ is equal to the order of $\t_B$. Therefore we see that
the rational elliptic surfaces of the types he considers, with
$n=2,3,4,6$ respectively, correspond to the cases in our list with
$n=2,3,4,6$, $m=n$ and $f_0 = I_0$. That is, respectively, cases
$\# 34, 35, 37$; $\# 30,31$; $\# 17$ and $\# 13$. The resulting Calabi-Yau threefolds correspond to the first line for each of these cyclic groups in table \ref{t:threefolds}.

\section{Outlook}\label{s:outlook}

In this paper we produced a complete classification of finite automorphism groups of rational elliptic surfaces such that they lift to free automorphism groups on the smooth fiber product of two rational elliptic surfaces. This work opens up various avenues of research.

\begin{itemize}
\item We classified and studied in detail a large class of finite automorphisms of rational elliptic surfaces. However, since our goal was ultimately to construct non-simply connected Calabi-Yau threefolds, we restricted ourselves to automorphisms acting freely on the fiber at infinity. An obvious direction of research would be to continue our study from a purely two-dimensional point of view and classify all finite automorphisms of rational elliptic surfaces.

\item As mentioned in the introduction, one of our simplifying assumptions was to restrict ourselves to the study of Calabi-Yau threefolds constructed as smooth fiber products of two rational elliptic surfaces. However, an important aspect of Schoen's work \cite{Schoen:1988} was precisely to study the more complicated situation when the fiber product has ordinary double point singularities. He showed that in many cases, the minimal resolution of the singular Calabi-Yau threefold is also Calabi-Yau. It would be interesting, both from a mathematical and a physics point of view, to construct free automorphism groups acting on these smooth resolutions, in order to construct an even larger family of non-simply connected Calabi-Yau threefolds.

\item From a physics point of view, the motivation behind this work is to construct a large class of non-simply connected Calabi-Yau threefolds suitable for string compactifications. In particular, we are interested in $E_8 \times E_8$ heterotic string vacua, which roughly speaking are given by triples $(X,V,G)$, where $X$ is a Calabi-Yau threefold and $V$ is a stable vector bundle on $X$ with structure group $G \subseteq E_8$. One way to obtain realistic four-dimensional standard model physics out of these compactifications is to consider non-simply connected Calabi-Yau threefolds $X$, with $G=SU(4)$ or $SU(5)$. In this work we produced a large class of non-simply connected $X$; the next step in this program is to construct and study ``standard-model bundles" on these threefolds (possibly by constructing invariant stable bundles on the simply connected Calabi-Yau covers $\tX$), by which we mean $SU(4)$ or $SU(5)$ bundles satisfying all the requirements needed to yield realistic four-dimensional physics. The only bundles known so far which give a consistent compactification
with the correct spectrum are those constructed in \cite{Bouchard:2005ag,Bouchard:2006dn} on the
Calabi-Yaus built from the $4$-dimensional family of rational elliptic
surfaces in case $\# 35$ of table \ref{t:list2}.
\end{itemize}

\appendix

\section{Sketch of the proof of remark \ref{r:as}}

In this appendix we sketch a case-by-case proof of remark \ref{r:as}, where
we claimed that the set of allowed sections $AS$ is non-empty for all lines
in tables \ref{t:list1} and \ref{t:list2}.

When the group is cyclic, the automorphism group is generated by an
automorphism of the form
\begin{equation}
\t_B = t_\xi \circ \alpha_B, \qquad \xi \in \ker(\PH_m),
\end{equation}
and $AS$ will be non-empty if we can find a section $\xi$ which intersects
$f_\infty$ at a torsion point of order $n = d m$. When the automorphism
group is non-cyclic, one of the generators
is of the form above, and the
other is translation by a torsion section $\eta$. In this case we must find
a section $\xi$ as above, and also make sure that all the sections
$\cp_i(\xi)$, $i=1, \ldots, n-1$ intersect $f_\infty$ at
torsion
points other than those where $\eta$ and its multiples do,
so that the
full automorphism group
acts freely on $f_\infty$.

The easiest cases are when the group is cyclic and $m=n$ is prime
($\IZ_5$, $\IZ_3$ and $\IZ_2$). In each of these cases we can take a
section $\xi$
corresponding to a point of
$\ker(\PH_m)$ which has
minimal length in the Mordell-Weil group.
Using the
height pairing, it is easy to show that $\xi$ is
disjoint from the zero section, hence must intersect $f_\infty$ at a
non-zero torsion point of order $n$. Therefore $AS$ is non-empty.

For the non-cyclic cases, we take a section $\xi$ as above, and check using
the height pairing that it is also disjoint from the torsion sections used
to construct the second cyclic subgroup. This works in all non-cyclic
cases, except case 10.

The $\IZ_6$ cases 14 and 15 and the $\IZ_4$ cases 19 and 20 can be treated
similarly. For the $\IZ_6$ case 13, we consider a minimal point $\xi_2
\in
\ker(\PH_2) \subset E_8$, and a minimal point
$\xi_3 \in \ker(\PH_3)  \subset E_8$. Both sections are disjoint from
the zero section, hence
$\xi_2$ must intersect
$f_\infty$ at a torsion point of order 2, while $\xi_3$
intersects $f_\infty$ at
a torsion point of order 3. So the sum $\xi_2 \boxplus \xi_3$ must
intersect $f_\infty$ at
a torsion point of order 6.

The remaining cases are 10, 17 and 18. For cases 17 and 18, we have $m=4$,
and we want to find a section $\xi$ intersecting $f_\infty$ at a torsion
point of order 4.
If we can find at least 4 sections $\xi_i \in \ker(\PH_4)$,
$i=1, \ldots, 4$, which are all disjoint from each other
and from the zero
section, that is $\xi_i \cdot \xi_j = - \delta_{ij}$, and $\xi_i \cdot
\sigma = 0$, for all $i,j$, then one of them must intersect $f_\infty$ at a
4-torsion point, since there are
only 3 non-zero torsion points of order 2 on
the smooth elliptic fiber $f_\infty$. Now for both cases 17 and 18, using
the height pairing and an
explicit basis for $E_7$ and $E_6^*$, we can find 4
minimal points satisfying these conditions.

Finally case 10. Here we have $m=4$ and want to construct a $\IZ_4 \times
\IZ_2$ automorphism group. Hence we want a section $\xi$ that intersects
$f_\infty$ at a
4-torsion point, and we also want to make sure that $\alpha_B \xi \boxplus
\xi$ does not intersect $f_\infty$ at the same point as the 2-torsion
section $\eta$. If
we can find at least 8 sections which are all disjoint from each other
and
from the zero section, then one of them, call it $\xi$, will necessarily
intersect $f_\infty$ at a 4-torsion point but  $\alpha_B \xi \boxplus
\xi$ will intersect $f_\infty$ at a 2-torsion point
different than the one where the 2-torsion section $\eta$ intersects
$f_\infty$. So let us try
to find 8 such sections.

Describe $D_4^*$ as the square lattice $\IZ^4$ plus the non-integral element
$\frac{1}{2}(1,1,1,1)$. The 24 minimal points of norm $1$ are $\pm$ the 4
unit vectors, and $\frac{1}{2}(\pm 1, \pm 1 ,\pm 1 ,\pm1 )$. Note that any
minimal point is disjoint from $\sigma$ and the torsion section $\eta$, so
we need to find 7 mutually disjoint sections projecting to minimal points of
the lattice. Note also that for case 10 there are 4 $I_2$ fibers, and that
any minimal section intersects two non-neutral components.

Consider first the three sections $\xi_1, \xi_2, \xi_3$ corresponding to
$(1,0,0,0)$, $\frac{1}{2}(1,1,1,1)$ and $\frac{1}{2}(1,1,1,-1)$, and the two
additional
sections $\xi_4 = \xi_2 \boxplus \eta$, $\xi_5 =\xi_3 \boxplus \eta$.
It is easy to show that these 5 sections are all disjoint,
since their mutual height pairing is always greater than 0. Now consider the
section $\xi'_6$ corresponding to $(0,1,0,0)$. It is clearly disjoint from
$\xi_i$ for $i=2,3,4,5$, and it is either disjoint from $\xi_1$ or
$\xi'_6 \cdot \xi_1 = 1$. But then
$\xi'_6 \boxplus \eta$ is disjoint from $\xi_1$. We can therefore choose
$\xi_6$ to be
either $\xi'_6$ or $\xi'_6 \boxplus \eta$, whichever is disjoint from
$\xi_1$.
We then know that $\xi_1$ and $\xi_6$ must intersect the
non-neutral components of the same two
$I_2$ fibers.
Finally, consider the section
$\xi'_7$ corresponding to $(0,0,1,0)$. Both it and $\xi'_7 \boxplus \eta$
are again clearly disjoint from
$\xi_i$ for $i=2,3,4,5$, so
by taking $\xi_7$ to be either $\xi'_7$ or $\xi'_7 \boxplus \eta$, it can
be assumed to be disjoint from $\xi_1$.
But then $\xi_7, \xi_1, \xi_6$ all  intersect the
non-neutral components of the same two
$I_2$ fibers. It follows that
$\xi_7$ is also disjoint from $\xi_6$.
Hence, we have constructed a set of 7 mutually disjoint sections
which are also
disjoint from $\eta$ and $\sigma$, demonstrating that $AS$ is non-empty
in this case as well.

\end{document}